\documentclass[12pt, leqno]{article}
\usepackage{amssymb}
\def\q \m#1#2{{\raise1pt\hbox{$#1$}\kern-1pt\big/
               \kern-1pt\raise-1pt\hbox{$#2$}}}

\def\bJ{{\rm \bf J}}
\def\bN{{\rm \bf N}}
\def\bR{{\rm \bf R}}
\def\bZ{{\rm \bf Z}}

\def\bC{{\rm {\bf C}}}
\def\ch{{\rm  ch}}
\def\bH{{\rm \bf H}}
\def\sR{{ \rm \scriptsize  \bf R}}
\def\sZ{{ \rm \scriptsize  \bf Z}}
\def\sC{{ \rm \scriptsize  \bf C}}
\def\sN{{ \rm \scriptsize  \bf N}}
\def\sm{{ \rm \scriptsize   mod}}

\newfam\msbfam
\font\twelmsb=msbm10 at 12pt
\font\tenmsb=msbm10 at 10 pt
\font\sevenmsb=msbm10 at 7pt

\textfont\msbfam=\twelmsb
\textfont\msbfam=\tenmsb
\scriptfont\msbfam=\sevenmsb

\newtheorem{thm}{Theorem}[section]
\newtheorem{lemma}{Lemma}[section]
\newtheorem{defn}{Definition}[section]
 
 \newtheorem{prop}{Proposition}[section]
 \newtheorem{cor}{Corollary}[section]

\newcommand{\C}{{\cal C}}

\newcommand{\norm}[1]{||{#1}||}

\begin{document}

\renewcommand{\theequation}{\thesection.\arabic{equation}}
\setcounter{equation}{0}

\centerline{\Large {\bf Rigidity and Vanishing Theorems in $K$-Theory I}}

\vskip 10mm  
\centerline{\bf Kefeng   LIU\footnote{Partially supported by the Sloan 
Fellowship and an NSF grant.}, 
Xiaonan MA\footnote{Partially supported by SFB 288.}
and Weiping ZHANG\footnote{Partially supported by NSFC, MOEC and the Qiu Shi Foundation.}}
\vskip 8mm

{\bf Abstract.} In this paper, we first establish a $K$-theory  version of the 
 equivariant family index theorem for a circle action, then use it to 
prove several rigidity and vanishing theorems on the equivariant $K$-theory level.\\

{\bf 0 Introduction}.
In [{\bf W}], Witten considered the indices of elliptic 
operators on the free loop space ${\cal L} M$ of a spin manifold $M$. 
In particular the index of  the formal signature operator on the loop space is 
exactly the elliptic genus of Landweber-Stong [{\bf LS}].
Motivated by physics, Witten made the conjectures about the 
rigidity of these elliptic operators which says that 
their $S^1$-equivariant indices on $M$ are independent of $g\in S^1$.  
See  [{\bf L}] for the history of the subject.

 These conjectures were proved by Taubes [{\bf T}], Bott-Taubes [{\bf BT}], 
and by Hirzebruch [{\bf H}] and Krichever [{\bf K}]. Many aspects of mathematics 
are involved in their proofs. Taubes used analysis of Fredholm operators
and Witten's interpretation of the Atiyah-Bott-Segal-Singer Lefschetz 
fixed point formula; Krichever used cobordism; Bott and Taubes and 
Hirzebruch used the Atiyah-Bott-Segal-Singer Lefschetz 
fixed point formula. 
In [{\bf Liu1}], it was observed that these rigidity theorems are 
consequences of their modular invariance. This allowed Liu ([{\bf Liu1, 2}])
to give a simple 
and unified proof, as well as various further generalizations,
 of the above conjectures of Witten. In particular, several new 
vanishing theorems were found in [{\bf Liu1, 2}]. 

In many situations in geometry, it is rather natural and necessary to 
generalize the above rigidity and vanishing theorems to the family case. 
For example, to use elliptic operators to study the fundamental groups of
 a manifold, one can reduce the question to a problem involving a family 
of elliptic operators. In [{\bf  LiuMa1}] and [{\bf  LiuMa2}], 
Liu and Ma proved several family rigidity and vanishing theorems. 
Such theorems contain higher level 
vanishing terms which will be useful for our understanding of the 
relationship  between group actions and fundamental groups.

 To be more precise, let $M,\ B$ be two compact smooth manifolds, 
and $\pi: M\to B$ be a  smooth fibration with compact fibre $X$.
Assume that a compact Lie group $G$ acts fiberwise on $M$, that is, the action
preserves each fiber of $\pi$. Let $P$ be a family of $G$-equivariant 
elliptic operators along the fiber $X$. Then the family index of $P$, 
in the sense of Atiyah and Singer [{\bf AS}], 
\begin{eqnarray} 
{\rm Ind} (P) = {\rm Ker } P - {\rm Coker } P \in K_{G} (B),
\end{eqnarray}
is well-defined.
Note that ${\rm Ind} (P)$ is a virtual $G$-representation.  
Let  $\widehat{G}$ denote the space of all complex irreducible 
representations of $G$. By [{\bf S}, Proposition 2.2], we have 
\begin{eqnarray}
{\rm Ind} (P) = \oplus_{V\in \widehat{G}} {\rm Hom}_G (V, {\rm Ind} (P))
\otimes V
\end{eqnarray}
with ${\rm Hom}_G (V, {\rm Ind} (P)) \in K(B)$.
We denote by $({\rm Ind} (P))^G\in K(B)$ the $G$-invariant part of 
${\rm Ind} (P)$.

A family of elliptic operator $P$ is said to be
{\em rigid on the equivariant Chern character level }
with respect to this $G$-action, if the equivariant Chern character 
$\ch_g ({\rm Ind} (P)) \in H^* (B)$ is independent of $g \in G$. 
If $\ch_g ({\rm Ind} (P))$ is identically zero for any $g$, then we say $P$ has 
{\em vanishing property on the equivariant Chern character level}. 
More generally, we say that $P$ is 
{\em rigid on the equivariant $K$-theory level}, if 
${\rm Ind} (P) = ({\rm Ind} (P))^G$. If this index is identically zero
in $K_{G}(B)$, then we say that $P$ has 
{\em vanishing property on the equivariant $K$-theory level}. 
To study rigidity and 
vanishing properties of ${\rm Ind} (P)$, it is clear that we  only need to restrict to 
the case where $G=S^1$. From now on we make the  assumption that $G=S^1$. 

Note that the $K$-theory level rigidity and vanishing properties are more 
subtle than those on the Chern character level. The reason is that, by taking
the Chern character, some torsion elements involved in the index bundle 
might be killed. Such torsion elements may appear in the study of fundamental 
groups, which we hope to pursue in the future.

In [{\bf  LiuMa1}], Liu and Ma proved 
that the elliptic  genera are actually rigid on the equivariant 
Chern character level. Several  vanishing theorems for  equivariant 
Chern characters of these index bundles are also proved in [{\bf LiuMa1}].
Motivated by the family rigidity theorem of [{\bf  LiuMa1}, Theorem 2.1], 
it is rather natural to expect that the   elliptic  genera have rigidity 
and vanishing properties on the equivariant $K$-theory level. 
The purpose of the present paper is to show that this is indeed the case.

To achieve our purpose, we first establish   a $K$-theory version 
of the   equivariant family index theorem 
 [{\bf LiuMa1}, Theorem 1.1] for $S^1$-action. 
However, we are not able to derive this formula directly
 by  applying the localization formula in the equivariant $K$-theory
as in [{\bf ASe}], as the localizing process will also kill the torsion element in
 $K_{S^1}(B)$. 
Here instead, we combine the analysis of Wu-Zhang [{\bf WuZ}, \S 3],
which in turn relies on the technique of Bismut-Lebeau [{\bf BL}],
 with a deformation 
trick of Zhang [{\bf Z}, \S 2], which allows us to avoid the small eigenvalues problem, to establish such a formula.

To prove the  main results of this paper, which are stated in Section 2.1, we will introduce some 
shift operators on certain  vector bundles over the fixed point set of the circle
action, 
and compare the index bundles after the shift operation. 
Then we get a recursive relation of these index bundles which will in turn lead us
to the final result.
This part is essentially a reformulation of 
 the basic ideas of Taubes [{\bf T}].
Our main observation here is that we can directly do  the shift operations on the 
fixed point set by applying the $K$-theory version of the equivariant
family index theorem. 
In this way, we avoid the construction of the Dirac operator on 
the normal bundle in the loop space, as well as the associated 
analysis on the Fredholm properties of these operators (cf. [{\bf T}]).
Consequently, some of the shifting operations we will construct are not 
the same as
that in [{\bf T}]. This simplifies the computation  significantly.
In fact, in some sense the proof we will present may be considered as a rather 
subtle $K$-theory version of the proof of Bott-Taubes [{\bf BT}].

In a subsequent paper we will use the method in this paper to 
prove several rigidity and vanishing theorems on the equivariant
$K$-theory level for Spin$^c$-manifolds,  and also 
for almost complex manifolds.

This paper is organized as follows. In Section 1, we prove a 
$K$-Theory version of the equivariant family index theorem for circle action. 
In Section 2, we prove  the rigidity and vanishing theorems of elliptic genera 
on the equivariant $K$-theory level. The proofs of the main results
in Section 2 are  
based on two intermediate results which will be proved in Sections 3 and 4
respectively.

Some of the results of this paper have been announced in [{\bf LiuMaZ}].

$\ $

{\bf Acknowledgements.} Part of this work was done while the authors were visiting the
Morningside Center of Mathematics in Beijing during the summer of 1999. 
The authors would like to 
thank the Morningside Center  for its hospitality. The second author would also like to 
thank the Nankai Institute of Mathematics for its hospitality.

\newpage

\section{ \normalsize  A $K$-theory version of the  equivariant family 
index theorem}
\setcounter{equation}{0}

In this section, we will prove a $K$-theory version 
of the equivariant family  index theorem
 [{\bf LiuMa1}, Theorem 1.1] for  $S^1$-action, which will play
a crucial role in the following sections. As an immediate consequence,
we obtain a $K$-theory version of  the famous 
$\widehat{A}$-vanishing theorem of Atiyah and Hirzebruch [{\bf AH}]
for compact connected spin manifolds with non-trivial $S^1$-actions.
As was pointed out in the introduction, different from the case of the usual index, 
we can not  get the $K$-theory index formula  directly
 by  applying the localization formula in  equivariant $K$-theory
as in [{\bf ASe}], as the localizing process will kill the torsion element in
 $K_{S^1}(B)$. So the formula we will derive is more precise and subtle for 
the $S^1$-action case.

This section is organized as follows: 
In Section 1.1, we state a $K$-theory version of the equivariant family
index theorem for 
$S^1$-action on a family of spin manifolds. 
 In Section 1.2, we prove Theorem 1.1 by applying the techniques of 
[{\bf BL}, Sections 8, 9], [{\bf WuZ}, Section 3] and a deformation trick 
in [{\bf Z}].  
In Section 1.3, we generalize Theorem 1.1 to somewhat more general situation. 
In particular, we obtain a $K$-theory version of  the equivariant family index
theorem for 
$S^1$-action on a family of Spin$^c$ Dirac operators.

\subsection{\normalsize A $K$-theory version of the equivariant family index theorem}

Let $M,\ B$ be two compact manifolds, and $\pi: M\to B$  a smooth fibration with 
compact fibre $X$ such that $\dim X = 2 l$. 
Let $TX$ denote the relative tangent bundle.
Let $W$ be a  complex vector bundle over $M$ and $h^W$  an Hermitian metric 
on $W$. 

Let $h^{TX}$ be a Riemannian metric on $TX$ and $\nabla^{TX}$ be the 
corresponding Levi-Civita  connection on $TX$ along the fibre $X$. 
Then the Clifford bundle 
$C(TX)$ is the bundle of Clifford algebras over $M$ whose fibre at 
$x\in M$ is the Clifford algebra $C(T_x X)$ of $(TX, h^{TX})$. 

We assume that the bundle $TX$ is spin over $M$.
 Let $S (TX)= S^+(TX) \oplus S^-(TX)$ be the spinor bundle of $TX$. 
We denote by $c(\cdot)$ the Clifford action of $C(TX)$ on $S(TX)$. 
Let $\{ e_i\}$ be an oriented orthonormal basis of $(TX, h^{TX})$,
let $\{ e^i\} $ be its dual basis. Let 
\begin{eqnarray}
\tau= i^l c(e_1) \cdots c(e_{2l})
\end{eqnarray}
be the involution of $S(TX)$. Then $\tau|_{S^\pm (TX)} = \pm 1$.

Let $\nabla^{S(TX) }$ be the Hermitian
connection on $S(TX)$ induced by $\nabla^{TX}$. 
Let $\nabla^W$ be  an Hermitian connection on $(W, h^W)$. 
Let $\nabla^{S(TX) \otimes W}$ be the connection on 
$S(TX) {\otimes} W$ along the fibre $X$:
\begin{eqnarray}
\nabla^{S(TX) \otimes W}= \nabla^{S(TX) } \otimes 1 + 1 \otimes \nabla^W.
\end{eqnarray}

For $b\in B$, we denote by $E_b, E_{\pm,b}$  the set of $\C^{\infty}$-sections 
of  $S(TX)\otimes W$, $S^{\pm}(TX) \otimes W$  over the fiber $X_b$.
We regard  $E_b$ as the fibre of a smooth ${\bf Z}_2$-graded infinite
 dimensional vector bundle over $B$. Smooth sections of $E$ over $B$ 
will be identified to smooth sections of $S(TX)\otimes W$ over $M$. 

\begin{defn} Define the twisted Dirac operator to be
\begin{eqnarray}\begin{array}{l}
D^X = \sum_ic(e_i) \nabla^{S(TX)\otimes W}_{e_i}.
\end{array}\end{eqnarray}
\end{defn}
Then $D^X$ is a family Dirac operator which acts fiberwise on the 
fibers of $\pi$. For $b\in B$, $D^X_b$ denote the restriction of $D^X$ to the 
fibre $E_b$. $D^X$ interchanges $E_+$ and $E_-$. Let $D^X_{\pm}$ be the  
restrictions of $D^X$ to $E_{\pm}$.

Now we assume that $S^1$  acts fiberwise  on $M$.
We will consider that $S^1$ acts as identity on $B$. Without loss of 
generality we can assume that $S^1$ acts on $(TX, h^{TX})$ isometrically.
 We also assume that the action of $S^1$ lifts to $S(TX)$ and $W$, and 
commutes with $\nabla^W$.

{}From [{\bf LiuMa1}, Proposition 1.1], we know that the difference bundle over $B$,
\begin{eqnarray}
{\rm Ind} (D^X) = {\rm Ker} D^X_{+,b}-  {\rm Ker} D^X_{-,b} ,
\end{eqnarray}
is well-defined in the equivariant $K$-group $K_{S^1}(B)$.

We denote $[n]$ ($n\in \bZ$)  the one dimensional complex vector space on which 
$S^1$ acts  as multiplication by $g^n$ for a generator  $g\in S^1$. 
By [{\bf S}, Proposition 2.2], we know that there exist $A\in \bN$, 
$R_n \in K(B)$  $(|n|\leq A)$ such that 
\begin{eqnarray}
{\rm Ind} (D^X) = \bigoplus_{|n|\leq A} R_n \otimes [n],
\end{eqnarray}
as an identification of $S^1$-difference bundles.

For $n\in \bZ$, let $E^n_{+,b},\ E^n_{-,b}$ be the subspaces of $E_{+,b},\ 
E_{-,b}$ where $S^1$ acts  as multiplication by $g^n$ for $g\in S^1$. 
Then we can consider $E^n_b= E^n_{+,b} \oplus  E^n_{-,b}$ as the fibres 
of a smooth $\bZ_2$-graded infinite dimensional vector bundle over $B$. 
Smooth sections of $E^n$ over $B$ will be identified to smooth sections of
 $S(TX) \otimes W$ over $M$ on which  
$S^1$ acts  as multiplication by $g^n$ for $g\in S^1$. 
Then the difference bundle over $B$,
\begin{eqnarray}
{\rm Ind} (D^X, n) = {\rm Ker} D^X_{+,b}|_{E^n_{+,b}}- 
 {\rm Ker} D^X_{-,b}|_{E^n_{-,b}},
\end{eqnarray}
is well-defined in the $K$-group $K(B)$. 

By the construction of the index bundle ${\rm Ind} (D^X)$ (cf.
 [{\bf LiuMa1}, Proposition 1.2]), one  knows well that
\begin{eqnarray}\begin{array}{ll}
{\rm Ind} (D^X, n)& =R_n,  \quad {\rm for}  \quad |n|\leq A,\\
& = 0, \quad  {\rm for}  \quad  |n| > A.
\end{array}\end{eqnarray}

Let $F=\{F_{\alpha}\} $ be the fixed point set of the circle action on $M$. 
Then $\pi: F_{\alpha}\to B$ (resp. $\pi: F\to B$) is a 
smooth fibration with fibre $Y_{\alpha}$ (resp. $Y$), and $Y$ is a 
totally geodesic compact submanifold in each fiber $X$. 
Let $\widetilde{\pi}: N\to F$ denote the normal bundle to $F$ in $M$. 
Then $N= TX/TY$. We identify $N$ as the orthogonal complement of $TY$ in $TX_{|F}$.  
Let $h^{TY},\ h^N$ be the corresponding metrics on $TY,\ N$ induced by $h^{TX}$.
Then, we have the following $S^1$-equivariant decomposition  of $TX$ when restricted to  $F$,
\begin{eqnarray}
TX_{|F} = N_{m_1} \oplus \cdots \oplus N_{m_l} \oplus TY, \nonumber
\end{eqnarray}
where each $N_{\gamma}$ is a complex vector bundle such that $g\in S^1$ acts on it 
by $g^{{\gamma} }$.  To simplify the notation, we will simply write that
\begin{eqnarray}\label{chern3}
TX_{|F} = \oplus_{v\neq 0} N_v \oplus TY, 
\end{eqnarray}
where $N_v$ is a complex vector bundle such that $g\in S^1$ acts on it 
by $g^v$ with  $v\in \bZ^*$.
Clearly, $N= \oplus_{v\neq 0} N_v$.
We will denote $N$ as a complex vector bundle, and $N_{\sR}$ the underlying 
real vector bundle, of $N$.

Since  $N_\sR$ is naturally oriented by the complex structure on $N$,
 $TY$ is naturally oriented by the orientations of $TX$ and $N$. Similarly, let
\begin{eqnarray}\label{chern4}
W_{|F} = \oplus_v W_{v}
\end{eqnarray}
be the $S^1$-equivariant decomposition of the restriction of $W$ over $F$. 
Here $W_v \ (v\in \bZ)$ is a complex vector bundle over $F$ on which 
$g\in S^1$ acts by $g^v$. 

Let $C(N_{\sR})$ be the Clifford algebra bundle of $(N_{\sR}, h^N)$. 
Then $\Lambda(\overline{N}^*)$ is a $C(N_{\sR})$-Clifford module. Namely, 
if $U\in N$, let $U' \in \overline{N}^*$ correspond to $U$ by 
the metric $h^{TX}$. If $U \in N$, we write 
\begin{eqnarray}\label{chern5}
 c(U)= \sqrt{2} U'\wedge, \quad c(\overline{U})=-\sqrt{2} i_{\overline{U}}.
\end{eqnarray}

On  $F$, let $\omega_2(TY),\ \omega_2(N)$ be the second 
Stiefel-Whitney classes of $TY,\ N$ respectively. Let $\det N_v$,
 $\det N= \otimes_v \det N_v$ be the determinant line bundles of $N_v,\ N$
over $F$. Let $c_1(\det N)$ be the first Chern 
class of $\det N$. As  $TX$ is spin, one gets
the following identity in $H^2(F, \bZ_2)$, 
\begin{eqnarray}\label{chern1}
\omega_2(TY)= \omega_2(N)= c_1(\det N) \ \ {\rm mod }\ (2).
\end{eqnarray}
As being explained in [{\bf LaM}, Appendix D, pp. 397], we can construct 
a Spin$^c$ -principal bundle and a complex spinor bundle 
$S(TY, (\det N)^{-1})$ over $F$ which  locally may be thought of as 
\begin{eqnarray}\begin{array}{l}
S\left(TY, (\det N)^{-1}\right)= S_0(TY) \otimes ( \det N)^{-1/2}, \\
\Lambda(\overline{N}^*)= S_0(N_{\sR})\otimes ( \det N)^{1/2},
\end{array}\nonumber
\end{eqnarray}
where $S_0(TY),\ S_0(N_{\sR})$ are the fundamental spinor bundles 
for the (possibly nonexistent) spin structure on $TY,\ N_{\sR}$, 
and where $( \det N)^{-1/2}$ is the 
(possibly nonexistent) square root of $(\det N)^{-1}$.

Since $Y$ is totally geodesic in $X$, the connection $\nabla^{TX}|_F$ 
also preserves the decomposition (1.8) of $TX$ over $F$.
Let $\nabla^{TY}$, $\nabla^N$, $ \nabla^{N_v}$
 be the corresponding induced connections on $TY$, $N$ and $N_v$,
Let $\nabla^{\Lambda(\overline{N}^*)}$ be the connection on 
$\Lambda(\overline{N}^*)$ induced by $\nabla^N$.

Let $\nabla^{(\det N)^{-1}}$ be the connection on $(\det N)^{-1}$ induced by 
$\nabla^N$. The connections $\nabla^{(\det N)^{-1}}$  and 
$\nabla^{TY}$ induce an  Hermitian connection  
$\nabla^{S(TY, (\det N)^{-1})}$ on $S(TY, (\det N)^{-1})$. 
In fact, locally, both $S_0(TY)$ and $ ( \det N)^{-1/2}$ exist 
and carry a canonical connection induced by 
$\nabla^{TY}$, $\nabla^{(\det N)^{-1}}$. We give the bundle $S(TY, (\det N)^{-1})$ 
the tensor product connection. It is standard that 
this connection is well defined globally (cf. [{\bf LaM}, Appendix D]).

Let $(e_1, \cdots, e_{2l'}),\ (e_{2l'+1}, \cdots, e_{2l})$ be the corresponding oriented 
orthonormal basis of $TY$ and $N_\sR$. Then $(e_1, \cdots, e_{2l})$ is an oriented 
orthonormal basis of $TX$. The $\bZ_2$-gradings on 
$S(TY, (\det N)^{-1})$, $\Lambda (\overline{N}^*)$ are defined by the involutions 
$i^{l'} c(e_1) \cdots c(e_{2l'})$, and $\tau^N=i^{l-l'} c(e_{2l'+1}) \cdots c(e_{2l})$ respectively. Also note that 
under this involution, 
$\Lambda^{\rm even} (\overline{N}^*)= (\Lambda (\overline{N}^*))^+$, 
$\Lambda^{\rm odd} (\overline{N}^*)= (\Lambda (\overline{N}^*))^-$.

{}From the  above discussion, we see that there is a natural isomorphism between 
${\bf Z}_2$-graded $C(TX)$-Clifford 
modules over $F$, 
\begin{eqnarray}\label{spinor1}
S\left(TY, (\det N\right)^{-1}) 
\widehat{\otimes} \Lambda\left(\overline{N}^*\right) \simeq S(TX)|_F.
\end{eqnarray}
Here we denote   the $\bZ_2$-graded tensor product by $\widehat{\otimes}$
(cf. [{\bf LaM}, pp. 11]). Furthermore, since $\nabla^N$ is $S^1$-invariant, 
one deduces easily that  
\begin{eqnarray}\label{spinor2}
\nabla^{S(TX)}|_F=\nabla^{S(TY, (\det N)^{-1})}\widehat{\otimes} 1 + 
 1 \widehat{\otimes} \nabla^{\Lambda(\overline{N}^*)}.
\end{eqnarray}

Let $V$ be an Hermitian vector bundle over $F$. Let $\nabla^{V}$ be 
an Hermitian connection on $V$. From now on, we will also write $D^Y \otimes V$ 
the family twisted Dirac operator on $S(TY, (\det N)^{-1}) \otimes V$ on $F$,
and $D^{Y_\alpha} \otimes V$ its restriction on $F_\alpha$.
Namely, let $\nabla^{S(TY, (\det N)^{-1}) \otimes V}$ be 
the tensor product connection on $S(TY, (\det N)^{-1}) \otimes V$ induced by 
$\nabla^{S(TY, (\det N)^{-1})}$ and  $\nabla^V$.
Then 
\begin{eqnarray}
D^{Y} \otimes V= 
\sum_{i=1}^{2l'} c(e_i) \nabla^{S(TY, (\det N)^{-1}) \otimes V}_{e_i}.
\end{eqnarray}

We use the notation 
\begin{eqnarray}\begin{array}{l}
{\rm Sym}_q (V) = \sum_{n=0}^{+\infty} q^n {\rm Sym}^n (V)\in K(F) [[q]],\\
\Lambda_q (V) = \sum_{n=0}^{+\infty} q^n \Lambda^n (V)\in K(F) [[q]] 
\end{array}\nonumber
\end{eqnarray}
for  the symmetric and exterior power   
operations in $K(F)[[q]]$ respectively.

Let us intorduce the notations: 
\begin{eqnarray}\begin{array}{l}
R(q) = q^{{1 \over 2} \sum_v |v| \dim N_v} \otimes_{v>0}
 \left( {\rm Sym}_{q^v} (N_v) \otimes \det N_v\right)\\
\hspace*{20mm} \otimes_{v<0}
  {\rm Sym}_{q^{-v}} \left(\overline{N}_v\right)
\otimes_v q^v W_v = \oplus _{n} R_n q^n,\\
R'(q) = q^{-{1 \over 2} \sum_v |v| \dim N_v} \otimes_{v>0}
  {\rm Sym}_{q^{-v}} \left(\overline{N}_v\right) \\
\hspace*{20mm} \otimes_{v<0}
 \left( {\rm Sym}_{q^{v}} (N_v)\otimes \det N_v\right)
\otimes_v q^v W_v = \oplus _{n} R'_n q^n.
\end{array}
\end{eqnarray}

Note that by [{\bf AH}], one knows that, as $TX$ is spin,
\begin{eqnarray}
\sum_v v \dim N_v =0 \quad  {\rm mod }\ (2).
\end{eqnarray}

We can now state the main result of this section as follows.

\begin{thm} For $n\in \bZ$, we have the following identity in $K(B)$,
\begin{eqnarray}   \label{spinor3}     \begin{array}{l}
{\rm Ind} \left(D^X, n\right) = \sum_\alpha (-1)^{\Sigma_{0<v} \dim N_v} 
{\rm Ind} \left(D^{Y_\alpha} \otimes R_n\right)\\
\hspace*{20mm}= \sum_\alpha (-1)^{\Sigma_{v<0} \dim N_v} 
{\rm Ind} \left(D^{Y_\alpha} \otimes R'_n\right).
\end{array}\end{eqnarray}
\end{thm}

If we take $W= \bC$ the trivial line bundle over $M$, then the operator 
$D^X \otimes \bC$ is exactly the canonical Dirac operator $D^X$. The following 
consequence generalizes the  famous 
$\widehat{A}$-vanishing theorem of Atiyah and Hirzebruch [{\bf AH}]
for compact connected spin manifolds with non-trivial $S^1$-actions
to the family case.
\begin{cor} If $M$ is connected, and the $S^1$ action is nontrivial, then
for the family of the canonical Dirac operators $D^X$ along the fibre $X$, 
one has
\begin{eqnarray}
{\rm Ind} \left(D^X\right) =0 \quad {\rm in}\quad  K_{S^1}(B).
\end{eqnarray}
\end{cor}

$Proof$: If the $S^1$ action is locally free, then by Theorem 1.1,
 we get directly (1.18). Otherwise, since on each $F_\alpha$,
\begin{eqnarray}
\sum_v |v| \dim N_v >0,
\end{eqnarray}
one deduces easily from Theorem 1.1 that 
\begin{eqnarray}
{\rm Ind} \left(D^X, n\right) =0 \ \quad {\rm in}\quad  K(B)
\end{eqnarray}
for any $ n\in {\bf Z} $,
from which (1.18) follows.\hfill $\blacksquare$\\

\subsection{\normalsize Proof of Theorem 1.1}

In this subsection, we  prove Theorem 1.1. The proof, which is
contained here for completeness of the present paper,
 is modeled on [{\bf WuZ}, Section 3] which in turn relies  on
the paper of Bismut and Lebeau [{\bf BL}]. 

This subsection is organized as follows. In Section 1.2.1, we  recall 
a result from  [{\bf WuZ}, Proposition 3.2] concerning the Witten 
deformation on flat 
space. In Section 1.2.2, we establish a Taylor expansion of 
$D^X + \sqrt{-1} Tc(H)$ near the fixed point set $F$, where $H$ is the
Killing vector field on $M$ generated by the circle action. In Section 1.2.3,
by using 
the techniques of [{\bf WuZ}, Section 3] and
[{\bf BL}, Section 9], we establish various estimates 
for certain operators induced from $D^X + \sqrt{-1}T c(H)$. In Section 1.2.4, 
we prove Theorem 1.1 by using a deformation trick in [{\bf Z}].

\subsubsection{\normalsize Witten's deformation on flat spaces}

Let $H$ be the canonical basis of ${\rm Lie} (S^1)= \bR$, i.e., for
$t\in \bR$, $\exp(tH)= e^{2 \pi i t} \in S^1$. Let $W$ be a complex
 vector space of dimension $n$ with an
Hermitian form. Let $\rho$ be a unitary representation of the circle 
group $S^1$ on $W$ such that all the weights are nonzero. 
Suppose $W^{\pm}$ are the subspaces of $W$ corresponding to positive and 
negative weights respectively, with $\dim_{\sC} W^- = \nu,\ \dim_{\sC} W^+ = n- \nu$. 
Let $z = \{z^i\}$ be the complex linear coordinates on $W$ such that 
the Hermitian structure on $W$ takes the standard form and 
$\rho $ is diagonal with weights $\lambda_i \in \bZ \setminus \{ 0 \}$
 $(1 \leq i \leq n)$, and $\lambda_i<0$ for $i\leq \nu$. 
The Lie algebra action is given by the vector field
 $H= 2 \pi \sqrt{-1} \Sigma_{i=1}^n 
\lambda_i (z^i {\partial \over \partial z^i}
 - \overline {z}^i {\partial \over \partial \overline{z}^i})$ on $W$. 
We write $K^{\pm}(W) = {\rm Sym} ((W^{\pm})^*) \otimes {\rm Sym}(W^\mp) 
\otimes \det(W^{\mp})$. 
Let $E$ be a finite dimensional complex vector space with an Hermitian 
form and suppose $E$ carries an unitary representation of $S^1$.

Let $\overline{\partial}$ be the twisted Dolbeault operator acting on
$\Omega^{0,*} (W, E) = \Gamma (\Lambda (\overline{W}^*)\otimes E)$, 
and $\overline{\partial}^*$ its formal adjoint. 
Let $D= \sqrt{2}(\overline{\partial} + \overline{\partial}^*)$. 
Let $c(H)$ be the Clifford action of $H$ on $\Lambda (\overline{W}^*)$ 
defined as in (1.10).  Let $L_H$ be the Lie derivative along $H$ acting on
$\Omega^{0,*} (W, E) $.

The following result was proved in [{\bf WuZ}, Proposition 3.2].

\begin{prop} 1. A basis of the space of $L^2$-solutions   of 
$D + \sqrt{-1} c(H)$ (resp. $D - \sqrt{-1} c(H)$) on the space of 
${\cal C}^\infty$ sections of $\Lambda (\overline{W}^*)$ is given by
\begin{eqnarray}
\Pi_{i=1}^\nu z_i^{k_i} \Pi^n_{i=\nu+1} \overline{z}_i^{k_i} 
e^{-\Sigma^n_{i=1}  \pi |\lambda_i| |z_i|^2} 
d \overline{z}_{\nu+1} \cdots d \overline{z}_n    \quad (k_i \in \bN)
\end{eqnarray}
with weight $\Sigma_{i=1}^\nu k_i |\lambda_i| + \Sigma_{i= \nu +1}^n 
(k_i+1)|\lambda_i|$ (resp.
\begin{eqnarray}
\Pi_{i=1}^\nu \overline{z}_i^{k_i} \Pi^n_{i=\nu+1} {z}_i^{k_i} 
e^{-\Sigma^n_{i=1}  \pi |\lambda_i| |z_i|^2} 
d \overline{z}_{1} \cdots d \overline{z}_\nu    \quad (k_i \in \bN)
\end{eqnarray}
with weight $-\Sigma_{i=\nu +1}^n k_i |\lambda_i|-\Sigma_{i= 1}^\nu 
(k_i+1)|\lambda_i|$).

So the space of $L^2$-solution of a given weight  of 
$D + \sqrt{-1} c(H)$ (resp. $D - \sqrt{-1} c(H)$) on the space of 
${\cal C}^\infty$ sections of $\Lambda (\overline{W}^*) \otimes E$ is finite 
dimensional. The direct sum of these weight spaces is isomorphic to 
$K^-(W) \otimes E$ (resp. $K^+(W)\otimes E$) as represnetations of $S^1$.

2. When restricted to an eigenspace of $L_H$, the operator $D+\sqrt{-1}c(H)$
(resp. $D - \sqrt{-1} c(H)$) has discrete eigenvalues.
\end{prop}

\subsubsection{\normalsize A Taylor expantion of certain deformed operators 
near the fixed-point set}

In this subsection, we will use the notation of Section 1.1.
We assume temporarily that $B$ is a point, and $Y$ is the fixed 
point set of the $S^1$ action on $X$.

Following [{\bf BL}, Section 8e)] and [{\bf WuZ}, Section 3.2], 
we now describe a coordinate 
system on $X$ near $Y$. 

For $\varepsilon >0$, set $B_\varepsilon = \{ Z\in N; |Z|< \varepsilon \}$. 
Since $X$ and $Y$ are compact, there exists $\varepsilon_0>0$ 
such that for $0 < \varepsilon \leq \varepsilon_0$, the exponential map 
$(y, Z) \in N \to
\exp_y^X (Z) \in X$ is a diffeomorphism  from $B_\varepsilon$ onto a 
tubular  neighborhood $U_\varepsilon$ of $Y$ in $X$.  
{}From now on, we identify $B_\varepsilon$ with $U_\varepsilon$ 
and use the notation $x= (y, Z)$ 
instead of $x = \exp _y ^X (Z)$. Finally, we identify $y\in Y$ 
with $(y, 0)\in N$.

Let $h^{TY}$, $h^N$ be the corresponding metrics on $TY$ and $N$ induced by $h^{TX}$. Let $dv_X,\ d v_{Y}$ and $d v_N$ be the corresponding volume elements
 on $(TX, h^{TX}),\ (TY, h^{TY})$ and $(N, h^N)$. 
Let $k(y,Z)\ ((y, Z)\in B_\varepsilon)$ be the smooth positive function  
defined by $dv_X(y,Z) = k(y,Z) dv_{Y}(y) dv_{N_y}(Z)$.
Then $k(y) = 1$ and ${\partial k \over \partial Z} (y) =0$ 
for $y\in Y$;
the latter follows from [{\bf BL}, Proposition 8.9] and the 
fact  that $Y$ is totally geodesic in $X$. 

For $x= (y, Z)\in U_{\varepsilon_0}$, we will identify $S(TX)_x$ 
with $S(TX)_y$ and $W_x$ with $W_y$ by the parallel transport with 
respect to the $S^1$-invariant connections $\nabla^{S(TX)}$ and $\nabla^W$
respectively, along the geodesic $t \to (y, tZ)$. The induced
identification of $(S(TX) \otimes W)_x$ with 
$(S(TX) \otimes W)_y$ preserves the metric and  the $\bZ_2$-grading, 
and is moreover $S^1$-equivariant. Consequently, $D^X$ can be considered as 
an operator acting on the sections of the bundle 
$\widetilde{\pi}^* (S(TX) \otimes W)|_{Y})$ 
over $B_{\varepsilon_0}$ commuting with the $S^1$ action.

For   $\varepsilon >0$, let $\bH (B_\varepsilon)$ 
(resp. $\bH(N)$) be the set of smooth sections of 
$\widetilde{\pi}^* (S(TX) \otimes W)|_Y$ on $B_\varepsilon$ (resp. 
on the total bundle of $N$). If $f,\ g\in \bH(N)$ have compact supports, 
we will write 

\begin{eqnarray}\label{inn1}
\left\langle f, g \right\rangle= \left(\frac{1}{2 \pi}\right)^{{\dim} X}
\int_{Y}   \left( \int_N \left \langle f,g \right \rangle (y,Z) 
dv_{N_y}(Z) \right)  dv _{Y}(y).
 \end{eqnarray}
Then $k^{1/2} D^X k^{-1/2}$ is a (formal) self-adjoint operator on 
$\bH(N)$. 

The connection $\nabla^N$ on $N$ induces a splitting $TN= N \oplus T^H N$, 
where $T^H N$ is the horizontal part of $TN$ with respect to $
\nabla^N$. Moreover, since $Y$ is totally geodesic, this splitting,
 when restricted to $Y$, is preserved by the connection 
$\nabla^{TX}$ on $TX|_Y$. Let $\widetilde{\nabla}$ be 
the connection on $(S(TX) \otimes W)|_Y$ induced by the 
restricion of $\nabla^{S(TX) \otimes W}$ to $Y$. 
We will still denote by $\widetilde{\nabla}$ the lift of the 
connection $\widetilde{\nabla}$ to
$\widetilde{\pi}^*(S(TX) \otimes W)|_Y)$,.

We choose a local orthonormal basis of $TX$ such that $e_1, \cdots, e_{2l'}$ 
form a basis of $TY$, and $e_{2l'+1}, \cdots, e_{2l}$, that of $N_\sR$. 
Denote the horizontal lift of $e_i$ $(1\leq i \leq 2l' )$ to $T^H N$ 
by $e_i ^H$. As in [{\bf BL}, Definition 8.16] and [{\bf WuZ}, (3.15)], we define
\begin{eqnarray}
D^H= \sum_{i=1}^{2l'} c(e_i) \widetilde{\nabla}_{e_i^H}, 
\qquad D^N= \sum_{i=2l'+1}^{2l} c(e_i) \widetilde{\nabla}_{e_i}.
\end{eqnarray}
Clearly, $D^N$ acts along the fibers  of $N$. Let $\overline{\partial}^N$ 
be the $\overline{\partial}$-operator  along the fibers of $N$, 
and let $\overline{\partial}^{N*} $ be its formal adjoint with respect to (1.23). 
By (1.13), it is easy to see that $D^N = \sqrt{2} (\overline{\partial}^N + 
\overline{\partial}^{N*}) $. Both $D^N$ and $D^H$ are self-adjoint
 with respect to (1.23).

For $T>0$, we define a scaling $f\in \bH(B_{\varepsilon_0}) \to 
S_T f \in \bH( B_{\varepsilon_0 \sqrt{T}})$ by
\begin{eqnarray}
S_T f(y, Z) = f \left(y, {Z \over \sqrt{T}}\right),  \qquad (y, Z) 
\in B_{\varepsilon_0 \sqrt{T}}.
\end{eqnarray}
For a first order differential operator 
\begin{eqnarray}
Q_T = \sum_{i=1}^{2l'} a_T^i (y,Z) \widetilde{\nabla}_{e_i^H} + 
\sum_{i=2l'+1}^{2l} b^i_T(y,Z) \widetilde{\nabla}_{e_i} + c_T(y, Z)
\end{eqnarray}
acting on $\bH(B_{\varepsilon_0 \sqrt{T}})$, where $a_T^i,\ b_T^i$ and $c_T$
 are endomorphisms of $\widetilde{\pi}^*((S(TX)\otimes W)|_Y)$, we write
\begin{eqnarray}
Q_T = O\left(|Z|^2 \partial ^N + |Z| \partial ^H + |Z| + |Z|^p\right), \qquad (p\in \bN),
\end{eqnarray}
if there is a constant $C>0$ such that for any $T\geq 1$, 
$(y,Z)\in B_{\varepsilon_0 \sqrt{T}}$, we have 
\begin{eqnarray}
\begin{array}{l}
\left|a_T^i (y,Z)\right| \leq C|Z| \qquad (1\leq i \leq 2l'),\\
\left|b_T^i (y,Z)\right| \leq C|Z|^2 \qquad ( 2l'+1 \leq i \leq 2l),\\
\left|c_T (y,Z)\right| \leq C\left(|Z|+ |Z|^p\right).
\end{array}
\end{eqnarray}

Let $J_H$ be the representation of Lie($S^1$) on $N$. Then $Z\to J_H Z$ 
is a Killing vector field on $N$. We have the following analog of
 [{\bf BL}, Theorem  8.18] and [{\bf WuZ}, Proposition 3.3].
\begin{prop} As $T\to + \infty$,
\begin{eqnarray}\begin{array}{l}
S_T k^{1/2} \left(D^X + \sqrt{-1} T c(H)\right) k^{-1/2} S^{-1}_T 
= \sqrt{T} \left(D^N + \sqrt{-1} c(J_H Z) \right) + D^H \\
\hspace*{30mm}
+ {1 \over \sqrt{T}} O\left(|Z|^2 \partial ^N + |Z| \partial ^H + |Z| + |Z|^3\right).
\end{array}
\end{eqnarray}
\end{prop}

{\em Proof} : Since $Y$ is totally geodesic in $X$
and the actions of $S^1$ on $N$ and $M$ commute with the exponential
map, one can proceed as in the proof of [{\bf WuZ}, Proposition 3.3] and
[{\bf BL}, Section 8] to
get (1.29).
\hfill $\blacksquare$\\

By Proposition 1.1, the solution space of the operator 
$D^N+ \sqrt{-1} c(J_HZ)$ along the fiber $N_y$ ($y\in F$) is 
(the $L^2$-completion of) $K^-(N_y)\otimes W_y$.
They form an infinite dimensional Hermitian holomorphic vector bundle
$K^-(N)\otimes W|_F$ over $Y$, with the Hermitian connection $\nabla^Y$
induced from those on $N$ and $ W|_Y\to Y$.

Let $\bH^0(Y)$ be the Hilbert space of square-integrable sections of
$S(TY, (\det N)^{-1})\otimes K^-(N)\otimes W|_Y$, and
$\bH^0(N)$, that of the bundle $\widetilde{\pi}^*((S(TX) \otimes W)|_Y)$, equipped with the
corresponding Hermitian forms.  By using (\ref{spinor1}), 
we define an embedding $\psi\ \colon\bH^0(Y)\to\bH^0(N)$ by
	\begin{eqnarray}\label{emb}
\psi\colon\alpha \otimes\beta\in\bH^0(Y) \to 
\widetilde{\pi}^*\alpha \wedge\tau(\beta)\in\bH^0(N).
	\end{eqnarray}
Here $\alpha\in \Gamma (Y, S(TY, (\det N)^{-1}))$, 
$\beta\in L^2(K^-(N)\otimes W|_Y)$
and $\tau$ is the isometry from $L^2(K^-(N)\otimes W|_Y)$ to
$L^2(\widetilde{\pi}^*(\Lambda (\overline{N}^*)\otimes W|_Y))$ given by
Proposition 1.1.
Clearly, $\psi$ is an isometry onto its image which we denote by $\bH'^{,0}$.
Let $p: \bH^0(N)\to\bH'^{,0}$ be the orthogonal projection.
Then we have the following analog of [{\bf BL}, Theorem 8.21]
and [{\bf WuZ}, Proposition 3.4], which can be proved in the same way as in 
[{\bf BL}] and [{\bf WuZ}].

\begin{prop} The following identity for operators acting on
${\bf H}^0(Y)$ holds,
	\begin{eqnarray}
\psi^{-1}p\,D^Hp\,\psi=D^Y \otimes R(1),
	\end{eqnarray}
where $R(1)$ is defined in (1.15).
\end{prop}

\subsubsection{ Estimates of the operators as 
$\mbox{\boldmath$\mathit T\to +\infty$}$}

We still assume temporarily that $B$ is a point. For $p\geq 0$, let $\bH^p(X)$, $\bH^p(N)$ and $\bH^p(Y)$ be 
the $p$-th Sobolev spaces of  sections of the bundles $S(TX)\otimes W \to X$,
$\widetilde{\pi}^* (S(TX)\otimes W )|_Y \to N$ and 
$S(TY, (\det N)^{-1}) \otimes K^-(N)\otimes W|_Y\to Y$ respectively.
The group $S^1$ acts on all these spaces. For any $\xi \in {\bf Z}$,
let $\bH^p_\xi(X)$,
$\bH^p_\xi(N)$ and $\bH^p_\xi(Y)$ be the corresponding subspaces
of weight $\xi\in \bZ$.
Recall that the constant $\epsilon_0>0$ is defined in last subsection.
We now take $\epsilon \in (0,\frac{\epsilon_0}{2}]$, which is small enough for each
eigenvalue of $L_H$ we will consider, but otherwise can be assumed to be fixed.
Let $\rho\colon\bR \to[0,1)$ be a smooth function such that
	\begin{eqnarray}
\rho(a)=\left \{ \begin{array}{l}
1 \quad {\rm if} \quad a \leq 1/2,\\
0 \quad {\rm if} \quad a \geq 3/4.

	\end{array} \right.
 \end{eqnarray}
For $Z\in N$, set $\rho_\epsilon(Z)=\rho(\frac{|Z|}{\epsilon})$.
Let $\alpha\in \Gamma(Y, S(TY, (\det N)^{-1}))$,
$\beta\in L^2(K^-(N)\otimes W|_Y)$.
We define a linear map $I_{T, \xi}:\bH_\xi^p(Y) \to {\bf H}_\xi^p(N) $ by
	\begin{eqnarray}
\alpha\otimes\beta\in\bH_\xi^p(Y) \longmapsto I_{T, \xi}\sigma=
\frac{\rho_\epsilon\norm{\beta}_0}{\norm{\rho_\epsilon S^{-1}_T
(\tau(\beta))}_0}
\,\widetilde{\pi}^*\alpha\wedge S^{-1}_T(\tau(\beta))\in {\bf H}_\xi^p(N) .
	\end{eqnarray}
Let the image of $I_{T,\xi}$ from $\bH^p_\xi(Y)$ be
$\bH^p_{T,\xi}(N)=I_{T,\xi}\bH^p(Y)\subset\bH^p_\xi(N)$.
Denote the orthogonal complement of $\bH^0_{T,\xi}(N)$ in $\bH^0_\xi(N)$
by $\bH^{0,\perp}_{T,\xi}(N)$,
and let $\bH^{p,\perp}_{T,\xi}(N)=\bH^p_\xi(N)\cap\bH^{0,\perp}_{T,\xi}(N)$.
Let $p_{T,\xi}$ and $p^\perp_{T,\xi}$ be the orthogonal projections from
$\bH^0_\xi(N)$ onto $\bH^0_{T,\xi}(N)$ and $\bH^{0,\perp}_{T,\xi}(N)$ respectively.

Since the bundle $S(TX)\otimes W$ over $U_{\epsilon_0}$ is 
identified with $\widetilde{\pi}^* (S(TY, (\det N)^{-1})\otimes  
\Lambda (\overline{N}^*) \otimes W|_Y)$ 
 over $B_{\epsilon_0}$, we can consider $k^{-{1/2}}I_{T, \xi}\sigma$
as an element of $\bH^p_\xi(X)$ for $\sigma\in\bH^p_\xi(Y)$.
Define the linear map $J_{T,\xi}$ by
	\begin{eqnarray}
\sigma\in\bH^p_\xi(Y)\longmapsto J_{T, \xi}\sigma=k^{-1/2}I_{T, \xi}\sigma\in\bH^p_\xi (X).
	\end{eqnarray}
Let $\bH^p_{T, \xi}(X)=J_{T, \xi}\bH^p(Y)$ be the image.
Denote the orthogonal complement of $\bH^0_{T, \xi}(X)$ in $\bH^0_\xi(X)$
by $\bH^{0,\perp}_{T, \xi}(X)$, 
and let $\bH^{p,\perp}_{T, \xi}(X)=\bH^p_\xi(X)\cap\bH^{0,\perp}_{T, \xi}(X)$.
Let $\bar{p}_{T,\xi}$ and $\bar{p}^\perp_{T,\xi}$ be the orthogonal projections from
$\bH^0_\xi(X)$ onto $\bH^0_{T, \xi}(X)$ and $\bH^{0,\perp}_{T, \xi}(X)$ respectively.
It is clear that $\bar{p}_{T, \xi}=k^{-1/2}p_{T, \xi} k^{1/2}$.

For any (possibly unbounded) operator $A$ on $\bH^0_\xi(X)$, write
	\begin{eqnarray}
A=\left(	\begin{array}{cc} 	A^{(1)} & A^{(2)} \\
					A^{(3)} & A^{(4)}
		\end{array}	\right)
	\end{eqnarray}
according to the decomposition
$\bH^0_\xi(X)=\bH^0_{T, \xi}(X)\oplus\bH^{0,\perp}_{T, \xi}(X)$, i.e.,
$A^{(1)}=\bar{p}_{T, \xi} A\,\bar{p}_{T, \xi}$,
$A^{(2)}=\bar{p}_{T, \xi} A\,\bar{p}^\perp_{T, \xi}$,
$A^{(3)}=\bar{p}^\perp_{T, \xi} A\,\bar{p}_{T, \xi}$ and
$A^{(4)}=\bar{p}^\perp_{T, \xi} A\,\bar{p}^\perp_{T, \xi}$.

Let $D_T=D^X+\sqrt{-1} Tc(H)$, where now $H$ denotes the Killing vector field
on $M$ generated by the circle action.
Let $D_{T, \xi}$ and $D^Y_\xi$ be the restrictions of the operators $D_T$ and
$D^Y\otimes R(1)$ on $\bH^0_\xi(X)$ and $\bH^0_\xi(Y)$ 
respectively. 

\begin{prop}\label{ESTIMATE}
1. As $T\to +\infty$,
	\begin{eqnarray}\label{A1}
J^{-1}_{T, \xi} D^{(1)}_{T, \xi} J_{T, \xi}=D^Y_\xi+O\left({1 \over \sqrt{T}}\right),
	\end{eqnarray}
where $O({1 \over \sqrt{T}})$ denotes a first order differential operator whose
coefficients are dominated by $\frac{C}{\sqrt{T}}$ ($C>0$).\\
2. For each $\xi\in\bZ$, there exists $C>0$ such that for any $T\geq 1$,
$\sigma \in\bH^{1,\perp}_{T, \xi}(X)$, $\sigma'\in\bH^1_{T, \xi}(X)$, we have
	\begin{eqnarray}\begin{array}{l}
\norm{D^{(2)}_{T, \xi}\sigma}_0
\leq C\left(\frac{\norm{\sigma}_1}{\sqrt{T}}+\norm{\sigma}_0\right),	
	\\
\norm{D^{(3)}_{T, \xi}\sigma'}_0
\leq C\left(\frac{\norm{\sigma'}_1}{\sqrt{T}}+\norm{\sigma'}_0\right).
	\end{array}\end{eqnarray}
3. For each $\xi\in\bZ$, there exist $\epsilon\in(0,\frac{\epsilon_0}{2}]$, 
$T_0>0$, $C>0$ such that for any $T\geq T_0$, 
$\sigma\in\bH^{1,\perp}_{T,\xi}(X)$, 
we have
	\begin{eqnarray}\label{A4}
\norm{D^{(4)}_{T, \xi}\sigma}_0\geq C\left(\norm{\sigma}_1+\sqrt{T}\norm{\sigma}_0\right).
	\end{eqnarray}
\end{prop}

{\em Proof} : Proposition 1.4 is the analogue of [{\bf WuZ}, Proposition 3.5]
and can be proved in the same way as in [{\bf WuZ}, pp. 165-166], which in turn 
relies on [{\bf BL}, Section 9].
\hfill $\blacksquare$

\subsubsection{\normalsize  Proof of Theorem 1.1}

We now go back to the family case. The important observation is that the 
analysis in the above two subsections works well to the fiberwise (twisted)
Dirac operators.

For any $u\in \bR$, we write 
\begin{eqnarray}
D_{T, \xi}(u) = D_{T, \xi}^{(1)} + D_{T, \xi}^{(4)} 
+ u\left(D_{T, \xi}^{(2)}+ D_{T, \xi}^{(3)} \right): E\to E.
\end{eqnarray}

The following lemma plays a key role in our proof of Theorem 1.1.

\begin{lemma} There exists $T_1>0$ such that for any $u \in [0,1]$ 
and $T\geq T_1$, $D_{T, \xi}(u)$ is a continuous family of Fredholm 
operators over $B$.
\end{lemma}

{\em Proof} : From Proposition 1.4, one deduces (cf.
[{\bf Z}, Lemma 2.2]) that there exist $C_1,\ C_2 >0$ such that 
for $u\in [0,1]$, $s\in E$ and $T$ large enough,
\begin{eqnarray}\label{ind1}
\norm{D_{T, \xi}s -D_{T, \xi}(u) s }_0\leq {C_1 \over \sqrt{T}} 
\norm{D_{T, \xi}s }_0 + C_2 \norm{s}_0.
\end{eqnarray}
{}From (\ref{ind1}) and the Fredholm property of $D_{T, \xi}$,
 one obtains the Fredholm property of $D_{T, \xi}(u)$ for sufficiently large $T$. 
\hfill $\blacksquare$\\

Recall that the index bundle construction [{\bf AS}] applies 
well to continuous families of Fredholm operators and that the 
homotopy invariance property for the index bundle still holds in this 
situation. 

{}From Lemma 1.1, one then gets the following identity of index bundles, when $T$
is large enough,
\begin{eqnarray}\label{ind2}
\begin{array}{l}
{\rm Ind}\left(D_{\xi}^X\right) = {\rm Ind}\left(D_{T, \xi}\right) 
= {\rm Ind}\left(D_{T, \xi}(0)\right)\\
\hspace*{10mm}   = {\rm Ind}\left(D_{T, \xi}^{(1)} \right) + {\rm Ind}
\left(D_{T, \xi}^{(4)} \right) 
\quad {\rm in }  \quad K(B),
\end{array}\end{eqnarray}
where in the last line,  ${\rm Ind}(D_{T, \xi}^{(1)} )$
(resp. $ {\rm Ind}(D_{T, \xi}^{(4)} ) $) is now regarded as a family
of Fredholm operators mapping from 
$\bH ^0_{T,\xi}(X)$ (resp. $\bH^{0,\perp}_{T, \xi}(X)$) to 
$\bH ^0_{T,\xi}(X)$ (resp. $\bH^{0,\perp}_{T, \xi}(X)$).

On the other hand, by the third part of Proposition 1.4, one has obviously that
\begin{eqnarray}\label{ind3}
{\rm Ind}\left(D_{T, \xi}^{(4)}\right)=0 \quad {\rm in }  \quad K(B),
\end{eqnarray}
when $T$ is large enough.

Let $D^{(1)}_{T, \xi, \alpha}$ be the restriction of $D^{(1)}_{T, \xi}$
on $F_\alpha$. 
{}From  (\ref{emb}), (\ref{A1}) and the definition of $J_{T,\xi}$, 
one deduces easily that when $T$ is large enough, one has,
\begin{eqnarray}\label{ind4}\begin{array}{l}
\displaystyle{
{\rm Ind} \left( D^{(1)}_{T, \xi}\right)=\sum_\alpha (-1)^{\Sigma_{0<v} \dim N_v} 
{\rm Ind} \left(J^{-1}_{T, \xi} D^{(1)}_{T, \xi, \alpha} J_{T, \xi}\right)  }\\
\hspace*{25mm}\displaystyle{
=\sum_\alpha (-1)^{\Sigma_{0<v} \dim N_v} 
{\rm Ind}  \left(D^{Y_\alpha}_\xi\right)
\quad {\rm in }  \quad K(B). }
\end{array}\end{eqnarray}
By (\ref{ind2}), (\ref{ind3}) and  (\ref{ind4}), 
one deduces the first equation of 
(\ref{spinor3}) easily.

To get the second equation of (\ref{spinor3}), we only need to apply 
the first equation of (\ref{spinor3}) to the $S^1$-action on $M$ 
defined by the inverse of the original $S^1$-action on $M$.

The proof of Theorem 1.1 is complete. 
\hfill $\blacksquare$\\

\subsection{ \normalsize The Spin$^c$ case}

We will keep the notations in Sections 1.1 and 1.2. For  future applications,
in this subsection, we will extend Theorem 1.1 to  Spin$^c$ cases.

Let $\pi: M\to B$ be a fibration of compact manifolds with compact fibre $X$ 
such that $\dim X = 2 l$
 and that $S^1$ acts fiberwise on $M$.  Let $h^{TX}$ be a metric on $TX$.
We assume that $TX$ is oriented.
Let $(W, h^W)$ be an Hermitian  complex vector bundle over $M$.

Let $V$ be a $2p$ dimensional oriented real vector bundle   over $M$. 
Let $L$ be a complex line bundle over $M$ with the property that the vector 
bundle $U= TX \oplus V$ obeys $w_2(U) = c_1(L) \ {\rm mod} \ (2)$ where $w_2$ 
denotes the second Stiefel-Whitney class, and $c_1(L)$ is the first Chern class of $L$. 
Then the vector 
bundle $U$ has a Spin$^c$-structure.  Let $h^V,\ h^L$ be metrics on $V,\ L$.
 Let $S(U,L)$ be the fundamental complex spinor bundle for $(U, L)$ 
[{\bf LaM}, Appendix D.9]. 

Assume that the $S^1$-action on $M$ lifts to $V$, $L$ 
and $W$, and assume the metrics $h^{TX},\ h^V,\ h^L,\ h^W$ are $S^1$-invariant. 
Also assume that  the $S^1$-actions on $TX,\ V,\ L$ lifts to $S(U,L)$.  

Let $\nabla^{TX}$ be the Levi-Civita connection on $(TX, h^{TX})$ along 
the fibre $X$. Let $\nabla^V$, $\nabla^L$ and $ \nabla^W$ be $S^1$ invariant and 
metric-compatible connections on $(V, h^V)$, $(L, h^L)$ and  $(W, h^W)$
respectively. 
Let $\nabla^{S(U, L)}$ be the Hermitian connection on $S(U, L)$ induced by 
$\nabla^{TX} \oplus \nabla^V$ and $\nabla^L$ as in Section 1.1. 
Let $\nabla^{S(U,L)\otimes W }$ be the 
tensor product connection on $S(U,L) \otimes W$ 
induced by $\nabla^{S(U,L)}$ and $\nabla^W$.

Let $\{e_i\}_{i=1}^{2l}$, $\{ f_j \}_{j=1}^{2p}$ be  the corresponding oriented orthonormal 
basis of $(TX, h^{TX})$ and $ (V, h^V)$.  Let $D^X$ be the family Spin$^c$-Dirac operator
 on the fiber $X$,
\begin{eqnarray}\label{dirac1}
D^X= \sum_{i=1}^{2l} c(e_i) \nabla^{S(U,L)\otimes W}_{e_i} .
\end{eqnarray}

There are two canonical  ways to consider $S(U,L)$ as a $\bZ_2$-graded vector 
bundle. Let 
\begin{eqnarray}\begin{array}{l}
\tau_s = i^{l } c(e_1) \cdots c(e_{2l}),\\
\tau_e = i^{l+p} c(e_1)\cdots c(e_{2l}) c(f_1) \cdots c(f_{2p})
\end{array}\end{eqnarray}
be two involutions of $S(U, L)$. Then $\tau_s^2= \tau_e^2=1$. We decompose 
$S(U,L)= S^+(U,L) \oplus S^-(U,L)$ corresponding to $\tau_s$ (resp. $\tau_e$) 
such that $\tau_s|_{S^{\pm}(U,L)} = \pm 1$ 
(resp. $\tau_e|_{S^{\pm}(U,L)} = \pm 1$). For $\tau = \tau_s $ or $\tau_e$, we 
 can define the index bundle 
${\rm Ind}_\tau (D^X) \in K_{S^1}(B)$ as in Section 1.1.

We have the following $S^1$-equivariant decomposition of $V$ restricted to $F$,
\begin{eqnarray}\label{dirac2}
V_{|F}= \oplus_{v\neq 0} V_v \oplus V_{0}^{\sR},
\end{eqnarray}
where $V_v$ is a complex vector bundle such that $ g$ acts on it by $g^v$, 
and $V_0^{\sR}$ is the real subbundle of $V$ such that $S^1$ acts as identity. For $v\neq 0$, let $V_{v,\bR}$ denote the underlying real 
vector bundle of $V_v$.
Denote  by  $2p'= \dim V_0^{\sR}$ and $ 2l'= \dim Y$. 

Let us write  
\begin{eqnarray} \label{dirac3}
L_F = L \otimes \left(\bigotimes _{v\neq 0} \det N_v \bigotimes _{v\neq 0} 
\det V_v
\right)^{-1}.
\end{eqnarray}
Then $TY \oplus V_0^{\sR}$ has a Spin$^c$ structure as
$w_2(TY\oplus V_0^{\sR}) = c_1(L_F) \ {\rm mod} \  (2)$. 
Let $S(TY \oplus V_0^{\sR}, L_F)$ be the fundamental spinor bundle for 
$(TY\oplus V_0^{\sR}, L_F)$.  

Let $D^Y,\ D^{Y_\alpha}$ be the families of Spin$^c$ Dirac operators  acting on
$ S(TY \oplus V_0^{\sR}, L_F)$ over $F,\ F_\alpha$. 
If $R$ is an Hermitian complex vector bundle equipped with an Hermitian
connection over $F$, let $D^Y\otimes R,\ D^{Y_\alpha}\otimes R$ denote 
the twisted Spin$^c$ Dirac operators on $ S(TY \oplus V_0^{\sR}, L_F) \otimes R$. 

Recall that $N_{v, \sR}$ and $V_{v,\sR}$  are canonically oriented by their
 complex structures. The decompositions (\ref{chern3}), ({\ref{dirac2}) 
induce  the orientations 
of $TY$ and $V_0^{\sR}$. Let $\{e_i\}_{i=1}^{2l'}$, $\{ f_j \}_{j=1}^{2p'}$ 
be the corresponding oriented orthonormal basis of $(TY, h^{TY})$ and $ (V_0^{\sR}, h^{V_0^{\sR}})$. 
There are two canonical ways to consider 
$ S(TY \oplus V_0^{\sR}, L_F) $ as a $\bZ_2$-graded vector bundle: let 
\begin{eqnarray}\begin{array}{l}
\tau_s = i^{l' } c(e_1) \cdots c(e_{2l'}),\\
\tau_e = i^{l'+p'} c(e_1)\cdots c(e_{2l'}) c(f_1) 
\cdots c(f_{2p'})
\end{array}\end{eqnarray}
be two involutions of $S(TY \oplus V_0^{\sR}, L_F)$. 
Then $\tau_s^2= \tau_e^2=1$. We decompose 
$S(TY \oplus V_0^{\sR}, L_F)= S^+(TY \oplus V_0^{\sR}, L_F)$ 
$\oplus S^-(TY \oplus V_0^{\sR}, L_F)$ corresponding to $\tau_s$ (resp. $\tau_e$) 
such that $\tau_s|_{S^{\pm}(TY \oplus V_0^{\sR}, L_F)} = \pm 1$ 
(resp. $\tau_e|_{S^{\pm}(TY \oplus V_0^{\sR}, L_F)} = \pm 1$).

By restricting to $F$, one has the following isomorphism of $\bZ_2$-graded 
Clifford modules over $F$,
\begin{eqnarray} \label{dirac4}
S(U, L) \simeq S\left(TY \oplus V_0^{\sR}, L_F\right) 
\widehat{\bigotimes _{v\neq 0} }\Lambda  N_v \widehat{\bigotimes _{v \neq 0}} \Lambda  V_v.
\end{eqnarray}
We denote by ${\rm Ind}_{\tau_s}$, ${\rm Ind} _{\tau_e}$ the index bundles 
corresponding to the involutions $\tau_s,\ \tau_e$.

Let $S^1$ act on $L$ by send $g\in S^1$ to $g^{l_c}$ $(l_c\in \bZ)$ on $F$. 
Then $l_c$ is locally constant on $F$. We define the following elements in $K(F) [[q^{1/2}]]$,
\begin{eqnarray}\quad \begin{array}{l} 
R_{\pm}(q) = q^{{1 \over 2}\Sigma_v |v| \dim N_v 
- {1 \over 2} \Sigma_v v \dim V_v +{1 \over 2} l_c}
 \otimes_{0<v}\left( {\rm Sym}_{q^v} (N_v) \otimes \det N_v\right)\\
\hspace*{20mm}\otimes_{v<0} {\rm Sym}_{q^{-v}} \left(\overline{N}_v\right) 
\otimes _{v\neq 0} 
\Lambda_{\pm q^v} (V_v) \otimes _v q^v W_v = \sum_n R_{\pm, n} q^n,\\
R'_{\pm}(q) = q^{-{1 \over 2}\Sigma_v |v| \dim N_v 
- {1 \over 2} \Sigma_v v \dim V_v +{1 \over 2} l_c} 
\otimes_{0<v} {\rm Sym}_{q^{-v}} \left(\overline{N}_v\right) \\
\hspace*{10mm}\otimes_{v<0}\left( {\rm Sym}_{q^{v}} ({N}_v) 
\otimes \det N_v \right)\otimes _{v\neq 0} \Lambda_{\pm q^v} (V_v) 
\otimes _v q^v W_v = \sum_n R'_{\pm, n} q^n.
\end{array}\end{eqnarray}

The following result generalizes [{\bf T}, Theorem 2.6] to the family case.

\begin{thm} For $n\in \bZ$, we have the following identity in $K(B)$,
\begin{eqnarray}\begin{array}{l} 
{\rm Ind}_{\tau_s} \left(D^X, n\right) = \sum_\alpha (-1)^{\Sigma_{0<v} \dim N_v} 
{\rm Ind} _{\tau_s} \left(D^{Y_\alpha} \otimes R_{+, n}\right)  \\
\hspace*{30mm}= \sum_\alpha (-1)^{\Sigma_{v<0} \dim N_v} 
{\rm Ind} _{\tau_s} \left(D^{Y_\alpha} \otimes R'_{+, n}\right) ,\\
 {\rm Ind}_{\tau_e} \left(D^X, n\right) = \sum_\alpha (-1)^{\Sigma_{0<v} \dim N_v} 
{\rm Ind} _{\tau_e} \left(D^{Y_\alpha} \otimes R_{-, n}\right)  \\
\hspace*{30mm}= \sum_\alpha (-1)^{\Sigma_{v<0} \dim N_v} 
{\rm Ind} _{\tau_e} \left(D^{Y_\alpha} \otimes R'_{-, n}\right).
\end{array}\end{eqnarray}
\end{thm}

{\em Proof }: The proof is a straightforward generalization of the proof 
of Theorem 1.1. The details are left to the interested reader. 
\hfill $\blacksquare$\\

\newpage

\section{ \normalsize Family rigidity and vanishing theorems}
\setcounter{equation}{0}

The purpose of this section is to establish the main results of this paper: 
the rigidity and vanishing theorems on the equivariant $K$-theory level.
The results in this section 
refine some of the results in [{\bf LiuMa1, 2}] to the $K$-theory level.

As in the previous sections, we let $\pi: M\to B$ be a  fibration of compact manifolds with fiber $X$. 
 We assume that  $S^1$ acts fiberwise on $M$ with fixed point set $F$,
and $TX$ has an $S^1$-equivariant spin structure. Then $\pi : F \to B$ is 
a fibration with fiber $Y$. 

Following Witten, we will introduce 
some elements $R(q) = \Sigma_{n\in \bN} q^n R_n \in K_{S^1} (M)[[q]]$.
 To prove the rigidity theorem for these elements, Taubes and 
Witten suggested 
to use some shift operators to get a relation like 
${\rm Ind} (D^X \otimes R_m, h) = {\rm Ind} (D^X \otimes R_{m+h}, h)$ 
for $h,\ m \in \bZ$. As $R_m=0$ for $m<0$, this implies the rigidity theorem.
See the paper of Taubes [{\bf T}] for a rigorous treatment.

To get a similar equality in the family case,  we first apply   
our $K$-theory version of the equivariant family index theorem 
[{\bf LiuMa1}, Theorem 1.1], Theorem 1.1, 
to reduce the problem to the fixed 
point set $F$. Then we introduce an auxiliary element in $K_{S^1}(F)[[q]]$. 
We study the  corresponding index bundles of the twisted Spin$^c$ Dirac 
operators on $Y$, which, after doing some shift operations, will be related to  
a term like ${\rm Ind} (D^X \otimes R_{m+h}, h)$. On the other hand, if we  apply
 our $K$-theory version of the equivariant family index theorem 
iterately, we may also relate the considered index bundle to 
a term like 
${\rm Ind} (D^X \otimes R_m, h) $. This then completes the  proof. To apply the equivariant family index theorem, we are inspired by
 the constructions 
of Taubes [{\bf T}, \S 6]. Namely, we will construct some operators on the fixed 
point set $M(n)$ of the induced $\bZ_n$-action on $M$, 
and apply the  equivariant family index theorem to them.

As was pointed out in the introduction, our main observation is that we 
can directly construct and apply the shift operators on the fixed point set. 
In this way, we avoid the use 
of the Dirac operators  on the normal bundle in the loop space 
and the associated  analysis of Fredholm properties of these operators 
 in [{\bf T}].

 This section is organized as follows: 
In Section 2.1, we state our main results, the rigidity and vanishing 
theorems on  the equivariant $K$-theory level.  In Section 2.2, we state two 
intermediate results which will be used to prove our main results stated in 
Section 2.1. 
In Section 2.3, we prove the family rigidity and vanishing theorem. 
In Section 2.4, we prove a vanishing theorem on the  equivariant $K$-theory 
level of the index bundle of the Dirac operator on loop spaces, which may be 
viewed  as a loop space analogue of  Corollary 1.1, and which also extends the
corresponding loop space analogue of the Atiyah-Hirzebruch 
 theorem in [{\bf Liu2}] to the family case.

Throughout this section, we use the  notations of Section 1.1.

\subsection{ \normalsize Family rigidity and vanishing theorems}

Let $\pi: M\to B$ be a  fibration of compact manifolds with fiber $X$ 
and $\dim X= 2l$. We assume that  $S^1$ acts fiberwise on $M$,
and $TX$ has an $S^1$-invariant spin structure. As in [{\bf AH}], 
by lifting to the double cover of $S^1$, we can assume that
 the second condition is always 
satisfied. Let $V$ be a real vector bundle  over $M$ with structure group 
Spin$(2k)$. We  assume that $V$ has an $S^1$-invariant spin 
structure.

The purpose of this section  is to prove that the elliptic operators 
introduced by 
Witten [{\bf W}] have some interesting rigidity and vanishing properties on the equivariant $K$-theory level. 
Let us recall some definitoins first.

For a complex (resp. real) vector bundle $E$ over $M$, let
\begin{eqnarray}\begin{array}{l}
{\rm Sym}_t (E) = 1 + t E + t^2 {\rm Sym}^2 E + \cdots ,\\
\Lambda_t (E) = 1 + tE + t^2 \Lambda^2 E + \cdots
\end{array}\end{eqnarray}
be the symmetric and exterior power   operations of $E$ 
(resp. $E\otimes_\sR {\bf C}$)
in $K(M)[[t]]$ respectively.
Set 
\begin{eqnarray}\begin{array}{l}
\Theta'_q(TX) = \otimes_{n=1}^\infty \Lambda_{q^n} (TX) 
\otimes _{n=1}^\infty {\rm Sym}_{q^n} (TX),\\
\Theta_q(TX) = \otimes_{n=1}^\infty \Lambda_{-q^{n-1/2}} (TX) 
\otimes _{n=1}^\infty {\rm Sym}_{q^n} (TX),\\
\Theta_{-q}(TX) = \otimes_{n=1}^\infty \Lambda_{q^{n-1/2}} (TX) 
\otimes _{n=1}^\infty {\rm Sym}_{q^n} (TX).
\end{array}\end{eqnarray}
We also define the following elements in $K(M) [[q^{1/2}]]$:
\begin{eqnarray}\begin{array}{l}
\Theta'_q(TX|V) = \otimes_{n=1}^\infty \Lambda_{q^n} (V) 
\otimes _{n=1}^\infty {\rm Sym}_{q^n} (TX),\\
\Theta_q(TX|V) = \otimes_{n=1}^\infty \Lambda_{-q^{n-1/2}} (V) 
\otimes _{n=1}^\infty {\rm Sym}_{q^n} (TX),\\
\Theta_{-q}(TX|V) = \otimes_{n=1}^\infty \Lambda_{q^{n-1/2}} (V) 
\otimes _{n=1}^\infty {\rm Sym}_{q^n} (TX),\\
\Theta^*_q(TX|V) = \otimes_{n=1}^\infty \Lambda_{-q^n} (V) 
\otimes _{n=1}^\infty {\rm Sym}_{q^n} (TX).
\end{array}\end{eqnarray}
Let $S(V) = S^+ (V) \oplus S ^- (V) $ be the spinor bundle of $V$.
 
Recall that the equivariant cohomology group $H^*_{S^1} (M, \bZ)$ 
of $M$ is defined by
\begin{eqnarray}
H^*_{S^1} (M, \bZ)= H^*(M \times_{S^1} ES^1, \bZ),
\end{eqnarray}
where $ES^1$ is the usual universal $S^1$-principal bundle over the 
classifying space $BS^1$  of $S^1$.
So $H^*_{S^1} (M, \bZ)$ is a module over $H^*(BS^1, \bZ)$ induced by the 
projection $\overline{\pi} : M\times _{S^1} ES^1\to BS^1$. 
Let $p_1(V)_{S^1},\ p_1(TX)_{S^1} \in H^*_{S^1} (M, \bZ)$ be the equivariant
 first Pontrjagin classes of $V$ and $TX$ respectively. 
As $V\times_{S^1}ES^1$ and $TX\times_{S^1}ES^1$ are spin over 
$M\times_{S^1}ES^1$, one knows that ${1\over 2}p_1(V)_{S^1}$ and
 ${1\over 2}p_1(TX)_{S^1}$ are well-defined in $H^*_{S^1}(M,\bZ)$
(cf. [{\bf T}, pp. 456-457]).
Also recall that 
\begin{eqnarray}\label{hyp0}
H^*(BS^1, \bZ)= \bZ [[u]]
\end{eqnarray}
with $u$ a generator of degree $2$.

In the following, we denote by $D^X\otimes W$ the family of Dirac operators
acting fiberwise  
on $S(TX) \otimes W$ as defined in Section 1.1. We also write 
$d_s^X=D^X\otimes S(TX)$. 

We can now state the main results of this paper as follows. The first one is a family 
generalization of the Witten rigidity theorems as being proved in [{\bf T}],
[{\bf BT}]  and [{\bf Liu1}]. It also refines  [ {\bf LiuMa1}, Theorem 2.1].

\begin{thm}(a) The family operators $d_s^X\otimes \Theta'_q(TX)$, 
$D^X \otimes \Theta_q(TX)$ and $D^X \otimes \Theta_{-q}(TX) $ are rigid
on the equivariant K-theory level.

(b) If  ${1\over 2}p_1(V)_{S^1} ={1\over 2} p_1(TX)_{S^1}$, then 
$D^X \otimes  ( S^+ (V) +  S ^- (V))\otimes\Theta'_q(TX|V) $,
 $D^X \otimes ( S^+ (V) -  S ^- (V))\otimes \Theta^*_q(TX|V)$, 
 $D^X \otimes \Theta_q(TX|V) $ and $D^X \otimes \Theta_{-q}(TX|V)$ are rigid
on the equivariant K-theory level.
\end{thm}

The second one generalizes the vanishing results of Taubes [{\bf T}, Proposition 10.1]
and Liu [{\bf Liu2}, Corollary 3.3] to the family case. It also refines a result in 
[{\bf LiuMa1}, Theorem 3.2].

\begin{thm}
 If  ${1\over 2}p_1(V)_{S^1} - {1\over 2}p_1(TX)_{S^1}= e \cdot  \overline{\pi}^* u^2$ 
with $e\in \bZ$ verifying $e<0$, then the index bundles of 
$D^X \otimes ( S^+ (V) + S ^- (V)) \otimes\Theta'_q(TX|V) $,
 $D^X \otimes ( S^+ (V) - S ^- (V))\otimes \Theta^*_q(TX|V)$, 
 $D^X \otimes \Theta_q(TX|V) $ and $D^X \otimes \Theta_{-q}(TX|V)$ are zero in
 $K_{S^1}(B)$. In particular, they are identically zero in $K(B)$.
\end{thm}

A quite interesting consequence of the above results is the following family 
$\widehat{A}$-vanishing theorem for loop spaces. It  extends the
corresponding loop space analogue of the Atiyah-Hirzebruch theorem
in [{\bf Liu2}, Theorem 6] to the family case.

\begin{thm} Assume $M$ is connected and the $S^1$-action is nontrivial. 
If ${1\over 2}p_1(TX)_{S^1} = -e \cdot \overline{\pi}^* u^2$ for some 
integer $e$, then the index bundle of 
$D^X\otimes_{n=1}^\infty {\rm Sym}_{q^n} (TX)$ 
as an element in $K_{S^1}(B)$, in particular,  as an element in $K(B)$, 
is identically zero. 
\end{thm}

As was pointed out by Dessai [{\bf De}], when the $S^1$-action is induced 
from a fiberwise  $S^3$ action on $M$ which preserves the Spin structure 
of $TX$, the condition
${1\over 2} p_1(TX)_{S^1} =- e \cdot \overline{\pi}^* u^2$ 
in $H^*_{S^1}(M,\bZ)$
   is  equivalent to 
${1\over 2} p_1(TX)=0$ in $H^*(M, \bZ)$. So 
one gets the following family vanishing theorem.

\begin{cor} Assume $M$ is connected and admits a 
nontrivial $S^1$ action induced by a fiberwise $S^3$-action 
which preserves the spin structure of $TX$. 
If ${1\over 2}p_1(TX) =0$, then
the index bundle as an element in $K_{S^1}(B)$ 
(in particular,  as an element in $K(B)$), of 
$D^X\otimes_{n=1}^\infty {\rm Sym}_{q^n} (TX)$, is identically zero. 
\end{cor}

Actually, our proof of these Theorems works under the following slightly 
weaker hypothesis. Let us first explain   some notations.

For each $n>1$, consider $\bZ_n \subset S^1$, the cyclic subgroup of order $n$.
We have the $\bZ_n$ equivariant cohomology of $M$ defined by 
$H^*_{\bZ_n}(M, \bZ) = H^*(M\times_{\bZ_n}  ES^1, \bZ)$, and there is 
a natural ``forgetful'' map $\alpha(S^1,\bZ_n): M\times_{\bZ_n} ES^1 
\to M \times_{S^1} ES^1$ which induces a pullback 
$\alpha(S^1, \bZ_n)^*: H^*_{S^1} (M, \bZ) \to H^*_{\bZ_n}(M, \bZ)$. 
The arrow which forgets the $S^1$ action altogether we denote by 
$\alpha(S^1, 1)$. Thus $\alpha(S^1, 1)^*: H^*_{S^1}(M, \bZ) \to H^*(M, \bZ)$
is induced by the inclusion of $M$ into $M\times_{S^1} ES^1$ as a fiber over 
$BS^1$. 

Finally, note that, if $\bZ_n$ acts trivially on a space $Y$, then there is 
another map $r^*: H^*(Y, \bZ) \to H^*_{\bZ_n} (Y, \bZ)$ induced by the 
projection $Y \times _{\bZ_n} ES^1 =Y \times B\bZ_n \stackrel{r}{\to} Y$.

We let $\bZ_{\infty}= S^1$. For each $1 < n \leq \infty$, let $i: M(n) \to M$
 be the inclusion of the fixed point set of $\bZ_n \subset S^1$ in $M$ 
and so $i$ induces $i_{S^1}: M(n) \times_{S^1} ES^1 \to M \times _{S^1} ES^1$.  

In  the rest of this paper, we suppose 
 that there exists some integer $e\in \bZ$ such that 
for each $1 < n \leq \infty$,
\begin{eqnarray}\label{hyp1}\begin{array}{l}
\alpha\left(S^1, \bZ_n\right)^* \circ i_{S^1} ^*\left({1\over 2} p_1(V - TX)_{S^1} 
- e \cdot  \overline{\pi}^* u^2\right) \\
\hspace*{10mm}=
 r^*\circ \alpha\left(S^1, 1\right)^*\circ i_{S^1} ^*\left({1\over 2} p_1(V-TX)_{S^1}
\right). 
\end{array}\end{eqnarray}

{\bf Remark 2.1.} The relation (\ref{hyp1}) clearly follows from the 
hypothesis  of Theorems 2.1 and 2.2 by pulling back and forgetting. 
Thus it is a weaker condition. 

$\ $

{\bf Remark 2.2.} If $e=0$, and $B$ is a point, 
(\ref{hyp1}) is exactly [{\bf BT}, (11.11)].\\

We can now state a slightly more general version of Theorems 2.1 and 2.2.

\begin{thm} Under the hypothesis (\ref{hyp1}), we have 

i) If $e=0$, then the index bundles of the elliptic operators in 
Theorem 2.1(b) are rigid on the equivariant K-theory level. 

ii) If $e<0$, then the index bundles of the elliptic operators in 
Theorem 2.1(b)  are zero as  elements in $K_{S^1}(B)$. In particular, these 
index bundles are zero in $K(B)$.
\end{thm}

The rest of this section is devoted to a proof of Theorem 2.4.

\subsection{\normalsize  Two intermediate results}

Let $F=\{F_ \alpha\}$ be the fixed point set of the circle action.
Then $\pi: F\to B$ is a fibration with compact fibre denoted by $Y = \{Y_ \alpha\}$.
Recall that 
\begin{eqnarray}
TX_{|F} = TY \oplus_{v\neq 0} N_{v, \sR}, \nonumber
\end{eqnarray}
where  $N_{v, \sR}$ denotes the underlying real bundle of the complex vector
 bundle $N_v$ on which $g\in S^1$ acts by multiplying by  $g^{v}$. 
 Since we can choose either $N_v$ 
or $\overline{N}_v$ as the complex vector bundle for $N_{v, \sR}$, in
what follows  
 we may and we will assume that
\begin{eqnarray}\label{hyp2}\begin{array}{l}
TX_{|F} = TY \oplus \bigoplus_{0<v} N_{v},\\
TX \otimes_\sR \bC = TY \otimes _\sR \bC 
\bigoplus_{0<v} (N_v \oplus \overline{N}_v),
\end{array}\end{eqnarray}
where $N_v$ is the complex vector bundle on which $g\in S^1$ acts by multiplying by $g^v$. Here some $N_v$ may be zero. Similarly, let
\begin{eqnarray}\label{hyp3}
V_{|F} = V_0^{\sR} \oplus \bigoplus_{0<v} V_v,
\end{eqnarray}
where $V_v$ is a complex  vector bundle on which $g\in S^1$ acts by multiplying by  
$g^v$, and $V_0^{\sR}$ is a real vector bundle on which $S^1$ acts as identity.

On $F$, let
\begin{eqnarray}\label{hyp4}\begin{array}{l}
e(N) = \sum_{0<v} v^2 \dim N_v, 
\qquad d' (N) =  \sum_{0<v} v \dim N_v,\\
e(V) = \sum_{0<v} v^2 \dim V_v, \qquad d' (V) =  \sum_{0<v} v \dim V_v. 
\end{array}\end{eqnarray}
Then $e(N),\ e(V),\ d'(N)$ and $ d'(V)$ are locally constant functions on $F$.

Let us introduce some line bundles:
\begin{eqnarray}\label{hyp5}\begin{array}{l}
L(N) = \otimes _{0<v} (\det N_v)^v,
\qquad  L(V) = \otimes _{0<v} (\det V_v)^v,\\
L= L(N)^{-1} \otimes L(V).
\end{array}\end{eqnarray}

We denote the Chern roots of $N_v$ by $\{ x_v ^j\}$ 
(resp. $V_v$ by $ u^j_v$),  and the Chern roots of 
$TY \otimes_{\sR} \bC$ by $\{ \pm y_j\}$ (resp. 
$V_0=V_0^{\sR} \otimes_{\sR} \bC$ by $\{ \pm  u_0^j\}$).
 Then if we take $\bZ_\infty = S^1$ in (\ref{hyp1}), we get
\begin{eqnarray}\label{hyp6}\begin{array}{l}
p_1(V)_{S^1}=\Sigma_{v,j}\left(u_v^j+vu\right)^2,\\
p_1(TX)_{S^1}=\Sigma_j(y_j)^2+\Sigma_{v,j}\left(x_v^j+vu\right)^2,\\
{1\over 2}\left(\Sigma_{v,j} \left(u_v^j + v u\right)^2
-  \Sigma_j (y_j)^2 - \Sigma_{v,j} \left(x_{v}^j + v u\right)^2 \right)- e u^2 \\
= {1\over 2}\left(\Sigma_{v,j} \left(u_v^j\right)^2 - 
\Sigma_j (y_j)^2 - \Sigma_{v,j} \left(x_{v}^j\right)^2\right) .
\end{array}\end{eqnarray}
By (\ref{hyp0}), (\ref{hyp6}), we get 
\begin{eqnarray}\label{hyp7}\begin{array}{l}
c_1(L) = \Sigma_{v,j} v u_{v}^j  -  \Sigma_{v,j}  v x_{v}^j =0,\\
e(V)- e(N) = \sum_{0<v} v^2 \dim V_v -\sum_{0<v} v^2 \dim N_v =2e.
\end{array}\end{eqnarray}
This means $L$ is a trivial complex line bundle over each component 
$F_\alpha$ of $F$, and $S^1$ acts on $L$ by multiplying by $g^{2e}$. 
So we can extend $L$ 
to a trivial complex line bundle over $M$, and we extend the $S^1$ action on it
 by multiplying the canonical section $1$ of $L$ to $g^{2e} \cdot 1$. 

The line bundles in (2.10) will play important roles in 
the next two sections which contain the proofs of Theorems 2.5, 2.6 to be stated below.

 We now define the following elements in $K_{S^1}(M) [[q^{1/2}]]$:
\begin{eqnarray}\begin{array}{l}
R_1(V)=\left(S^+ (V)+S^- (V)\right) \otimes_{n=1}^\infty \Lambda_{q^n} (V) ,\\
R_2(V)=\left(S^+ (V)- S^- (V)\right)\otimes_{n=1}^\infty \Lambda_{-q^n} (V), \\
R_3(V)=\otimes_{n=1}^\infty \Lambda_{-q^{n-1/2}} (V), \\
R_4(V)=\otimes_{n=1}^\infty \Lambda_{q^{n-1/2}} (V).
\end{array}\end{eqnarray}

In what follows, if $R(q)=\sum_{m \in {1 \over 2} \bZ} R_m q^m \in K_{S^1}(M) 
[[q^{1/2}]]$, 
we will also denote  ${\rm Ind} (D^X \otimes R_m, h)$ by
${\rm Ind} (D^X \otimes R(q), m, h)$. 

We first state a result which expresses the global equivariant family index via
the family indices on the fixed point set.

\begin{prop} For $m\in {1\over 2} \bZ$, $h\in \bZ$, $1\leq i \leq 4$,
 we have the following identity in $K(B)$,
\begin{eqnarray}
\begin{array}{l}
{\rm Ind} \left(D^X \otimes _{n=1}^\infty {\rm Sym}_{q^n} (TX)
 \otimes R_i(V), m, h\right)\\
= \sum_\alpha (-1)^{\Sigma _{0<v} \dim N_v}{\rm Ind} 
(D^{Y_\alpha} \otimes _{n=1}^\infty {\rm Sym}_{q^n} (TX) \otimes R_i(V)
\\
\hspace*{45mm}
 \otimes {\rm Sym} (\oplus_{0<v} N_v) \otimes _{0<v} \det N_v, m, h)
\end{array}
\end{eqnarray}
\end{prop}

{\em Proof.}   Proposition 2.1 follows directly from Theorem 1.1. \hfill $\blacksquare$\\

For $p\in \bN$, we define the following elements in $K_{S^1}(F) [[q]]$:
\begin{eqnarray} \label{hyp8}\begin{array}{l}
{\cal F}_p (X) = \bigotimes_{0<v}\Big (  \otimes_{n=1}^{\infty} 
{\rm Sym}_{q^n} (N_v) \otimes_{n> pv} {\rm Sym}_{q^n} (\overline{N}_v)  \Big)
\otimes_{n=1}^{\infty} {\rm Sym}_{q^n} (TY),\\
{\cal F}'_p (X) = \bigotimes_{\stackrel{0<v}{0\leq n \leq pv}}
\Big (  {\rm Sym}_{q^{-n}} (N_v)  \otimes \det N_v \Big ),\\
\\
{\cal F}^{-p} (X) = {\cal F}_p (X) \otimes {\cal F}'_p (X).
\end{array}\end{eqnarray}
Then, by (\ref{hyp2}),  over $F$, we have 
\begin{eqnarray}
{\cal F}^0 (X) = \otimes _{n=1}^\infty {\rm Sym}_{q^n} (TX)\otimes {\rm Sym} (\oplus_{0<v} N_v) \otimes _{0<v} \det N_v.
\end{eqnarray}

We now state two intermediate results on the relations between the family indices on the fixed point set. They will be used in the next subsection to 
prove Theorem 2.4.

\begin{thm} For $1\leq i \leq 4$, $h,\ p \in \bZ$ and $p>0$, 
$m\in {1\over 2} \bZ$,  we have 
\begin{eqnarray}\begin{array}{l}
\sum_\alpha (-1)^{\Sigma_{0<v} \dim N_v} {\rm Ind} (D^{Y_\alpha} 
\otimes {\cal F}^0 (X)  \otimes R_i(V), m , h)\\
=\sum_\alpha (-1)^{\Sigma_{0<v} \dim N_v} {\rm Ind} (D^{Y_\alpha} 
\otimes {\cal F}^{-p} (X)  \otimes R_i(V),\\
\hspace*{40mm} m+ {1\over 2} p^2 e(N) +{1\over 2} p d'(N), h).
 \end{array} \end{eqnarray}
\end{thm}

\begin{thm} For each $\alpha$, $1\leq i \leq 4$, $h,\ p \in \bZ$, $p>0$, 
$m\in {1\over 2} \bZ$, 
 we have the following identity in $K(B)$,
\begin{eqnarray}\begin{array}{l}
{\rm Ind} (D^{Y_\alpha} \otimes {\cal F}^{-p} (X)  \otimes R_i(V),
m + {1\over 2} p^2 e(N) +{1\over 2} p d'(N), h)\\
={\rm Ind} (D^{Y_\alpha} \otimes {\cal F}^0 (X)  \otimes R_i(V),
m + ph + p^2 e, h + 2pe).
\end{array} \end{eqnarray}
\end{thm}

Theorem 2.5 is a direct consequence of Theorem 2.7 to be stated below,
which will be proved in Section 4, while Theorem 2.6 will be proved in
Section 3.

To state Theorem 2.7,
set $J=\{v \in \bN|$  there exists $\alpha$ such that $N_v\neq 0$ on  $F_\alpha\}$ and 
\begin{eqnarray}\label{e1}
\Phi = \{ \beta \in ]0, 1]| {\rm there \  
exists } \ v\in J \ {\rm such\ that} \ \beta v\in \bZ \}.
\end{eqnarray}
We order the elements in  $\Phi $ so that
$\Phi=\{ \beta_i| 1\leq i \leq J_0, J_0\in \bN 
\ {\rm and} \ \beta_i < \beta_{i+1}\}$.
Then for any integer $1\leq i\leq J_0$,
 there exist $p_i,\ n_i \in \bN,\ 0< p_i\leq n_i$ with $ (p_i, n_i)=1$ such that 
\begin{eqnarray}\label{e2}
\beta_i = {p_i/n_i}. 
\end{eqnarray} 
Clearly, $\beta_{J_0}=1$. We also set  $p_0=0$ and $\beta_0=0$.

For $ 1\leq j \leq J_0$, $p \in \bN^*$, set
\begin{eqnarray}\label{e3}\begin{array}{l}
I^p_0 = \phi, \mbox{the empty set},\\
\displaystyle{I^p_j= \{ (v,n)| v\in J, (p-1)v<n \leq  pv,   
{n \over v}= p-1 + {p_j \over n_j} \} ,}\\
\displaystyle{\overline{I}^p_j= \{ (v,n)| v\in J, (p-1)v<n \leq pv, 
{n \over v} > p-1 +  {p_j \over n_j}\}.}
\end{array}\end{eqnarray}
For $0\leq j \leq J_0$, we write  
\begin{eqnarray}\label{e4}
\\
{\cal F}_{p,j} (X) = {\cal F}_p(X)\otimes {\cal F}'_{p-1} (X) 
\bigotimes_{(v,n)\in \cup_{i=1}^{j} I_i^p} \Big({\rm Sym}_{q^{-n}} ( N_v) 
\otimes  \det N_v \Big )\bigotimes_{(v,n)\in \overline{I}^p_{j}} 
{\rm Sym}_{q^{n}} ( \overline{N}_v) .\nonumber
\end{eqnarray}
Then 
\begin{eqnarray}\label{e5}\begin{array}{l}
{\cal F}_{p,0}(X) = {\cal F}^{-p+1}(X), \\
{\cal F}_{p,J_0}(X) = {\cal F}^{-p}(X).
\end{array}\end{eqnarray}

For $s\in \bR$, let $[s]$ denote the greatest integer which is less than or 
equal to the given number $s$. Set
\begin{eqnarray}\label{e6} \begin{array}{l}
e(p, \beta_j, N) = {1\over 2} \sum_{0<v} (\dim N_v )
\Big ((p-1) v + [{p_j v \over n_j}]\Big ) 
\Big ((p-1) v + [{p_j v \over n_j}]+1\Big ),\\
d'(\beta_j, N) =  \sum_{0<v} (\dim N_v )[{p_j v \over n_j}].
\end{array} \end{eqnarray}
Then $e(p, \beta_j, N), d'(\beta_j, N)$ are locally constant functions on $F$. 
And 
\begin{eqnarray}\label{e7}\begin{array}{l}
e(p, \beta_0, N)={1 \over 2} (p-1)^2  e(N) + {1 \over 2} (p-1) d'(N),\\
 e(p, \beta_{J_0}, N)={1 \over 2} p^2  e(N) + {1 \over 2} p d'(N).
\end{array} \end{eqnarray}

\begin{thm} For  $1\leq i \leq 4$, $1\leq j \leq J_0$, $p \in \bN^*$, 
$ h\in \bZ$, $m \in {1 \over 2} \bZ$, we have 
\begin{eqnarray}\label{e8}\begin{array}{l}
\sum_\alpha (-1)^{d'(\beta_{j-1}, N)  + \Sigma_{0<v} \dim N_v} 
{\rm Ind} (D^{Y_\alpha} \otimes {\cal F}_{p,j-1} (X) 
 \otimes R_i(V),\\
\hspace*{50mm}m +e(p, \beta_{j-1}, N), h)\\
=\sum_\alpha (-1)^{d'(\beta_{j}, N)  + \Sigma_{0<v} \dim N_v} 
{\rm Ind} (D^{Y_\alpha} \otimes {\cal F}_{p, j} (X)
  \otimes R_i(V),\\
\hspace*{50mm}m + e(p, \beta_j, N) , h).
\end{array}\end{eqnarray}
\end{thm}

{\em Proof}: The proof is given in Section 4.\hfill $\blacksquare$\\

{\em Proof of Theorem 2.5} : From  (\ref{e5}), (\ref{e7}), and Theorem 2.7, 
for $1\leq i \leq 4$,  $h\in \bZ$, $p\in \bN^*$ and  $m\in  {1 \over 2} \bZ$, 
we get the following identity in $K(B)$:
\begin{eqnarray}
\begin{array}{l}
\sum_\alpha (-1)^{ \Sigma_{0<v} (v+1) \dim N_v} {\rm Ind} 
(D^{Y_\alpha} \otimes {\cal F}^{-p} (X)  \otimes R_i(V),\\
\hspace*{40mm} m +{1\over 2} p^2 e(N) +{1\over 2} p d'(N), h)\\
=\sum_\alpha (-1)^{\Sigma_{0<v}  \dim N_v}
{\rm Ind} (D^{Y_\alpha} \otimes {\cal F}^{-p+1} (X) \otimes 
R_i(V), \\
\hspace*{40mm} m + {1\over 2} (p-1)^2 e(N) +{1\over 2} (p-1) d'(N), h).
\end{array}
\end{eqnarray}
As $\sum_{0<v} v \dim N_v =0 \ \sm \ (2)$, 
we get Theorem 2.5. \hfill $\blacksquare$

\subsection{ \normalsize Proof of Theorem 2.4}

{}From  Proposition 2.1, Theorems 2.5 and  2.6, for $1\leq i \leq 4$, 
$h, p\in \bZ$, $p> 0$ and $m\in {1 \over 2} \bZ$,  we get the following identity 
in $K(B)$,
\begin{eqnarray}\label{e9}\begin{array}{l}
{\rm Ind} (D^X \otimes _{n=1}^\infty {\rm Sym}_{q^n} (TX) \otimes R_i(V),
 m, h)\\
\hspace*{10mm}
={\rm Ind} (D^X \otimes _{n=1}^\infty {\rm Sym}_{q^n} (TX) \otimes R_i(V), 
m', h+ 2pe).
\end{array}\end{eqnarray}
with 
\begin{eqnarray}\label{e10}
m'=m+ ph+ p^2 e.
\end{eqnarray}

Note that, from (2.1) and (2.13), if $m<0$, for $h \in \bZ$, we have   
\begin{eqnarray}\label{e11}
{\rm Ind} (D^X \otimes _{n=1}^\infty {\rm Sym}_{q^n} (TX) \otimes R_i(V),
 m, h) =0 \quad {\rm in } \quad  K(B).
\end{eqnarray}

i) Assume that $e=0$. Let $h\in \bZ,\ m_0\in {1 \over 2} \bZ$ and  $h\neq 0$ be 
fixed. If $h>0$, 
we take $m'=m_0$, then for $p$ big enough, we get $m<0$ in (\ref{e10}). If $h<0$,  we take $m=m_0$,  then  for $p$ big enough, 
we get $m'<0$ in (\ref{e10}). From (\ref{e11}), we know that for $h\neq 0$, 
$m_0\in {1 \over 2} \bZ$ and $1\leq i \leq 4$, 
we get 
\begin{eqnarray}\label{e12}
{\rm Ind} (D^X \otimes _{n=1}^\infty {\rm Sym}_{q^n} (TX) \otimes R_i(V), m_0, h) 
=0 \quad {\rm in } \quad  K(B).
\end{eqnarray}

ii) Assume that  $e <0$.  For $h\in \bZ$, $m_0\in {1 \over 2} \bZ$, 
$1\leq i \leq 4$, we take $m=m_0$, then  for $p$ big enough, 
we get $m'<0$ in (\ref{e10}). From (\ref{e11}), we  again get (\ref{e12}).

The proof of Theorem 2.4 is complete. \hfill $\blacksquare$

$\ $

{\bf Remark 2.3.} It might be suitable to add a remark on the comparison of the various 
proofs of the Witten rigidity theorem given in [{\bf T}], [{\bf BT}], [{\bf Liu1}] 
and  the present paper. On one hand, the proofs in 
[{\bf BT}] and [{\bf Liu1}] relies on the Atiyah-Bott-Segal-Singer fixed point formula
(cf. [{\bf ASe}], [{\bf ASi}]) and the elliptic function theory, so they don't work on 
the $K$-theory level. This is reflected in [{\bf LiuMa1}] where a proof of a family Witten
rigidity theorem on the equivariant Chern character level is given by extending
the method in [{\bf Liu1}]. On the other hand,
consider the proofs given in [{\bf T}], [{\bf BT}] and  the present paper
(for the last see  Section 4 for more details). All these proofs  
rely on Taubes's idea that in certain steps one needs to consider Dirac operators
on the fixed point set of the induced ${\bf Z}_n$-actions. This requires that the topological
conditions imposed  in these proofs are for the half of the equivariant first Pontryagin
classes. While the proof in [{\bf Liu1}] works directly on the fixed point set of
the $S^1$-action, and thus works under weaker conditions on the equivariant first
Pontryagin classes without the factor ${1\over 2}$. This leads to the natural 
question  that whether there would exist a $K$-theory version of the proof in [{\bf Liu1}]. Very likely one needs to use Hecke operator in the theory of modular forms to understand the shift operators and the modular transformations.

\subsection{ \normalsize Proof of Theorem 2.3}
In fact, by setting $V=0$ in (\ref{hyp7}), we know that 
\begin{eqnarray}
\Sigma_{0<v} v^2\dim N_v =-2e.
 \end{eqnarray}
Thus  the case $e>0$ can never happen. If $e=0$, then all the numbers 
$\dim N_v $
are zero, so that the $S^1$-action cannot  have  fixed points. 
{}From Theorem 1.1, we know that the index bundle is zero in $K_{S^1}(B)$. 
For $e<0$, one may take $V=0$ in Theorem 2.4 to derive  Theorem 2.3.

The proof of Theorem 2.3 is complete. \hfill $\blacksquare$\\

\newpage

\section{ \normalsize Proof of Theorem 2.6}
\setcounter{equation}{0}

In this section, we will prove Theorem 2.6 
by  introducing some shift operators as in [{\bf T}, \S 7].

This section is organized as follows: 
In Section 3.1, we introduce some notations. In Section 3.2, we prove 
Theorem 2.6 by  introducing some shift operators as in [{\bf T}, \S 7].

Throughout this section, we use the  notations of Section 2.

\subsection{ \normalsize Reformulation  of Theorem 2.6}

To simplify the notations, we  introduce some new notations in this subsection.

For $n_0\in \bN^*$, we define a number operator $P$ on 
$K_{S^1}(M)[[q^{1 \over n_0}]]$ in the following way: if 
$R(q)= \oplus_{n\in {1 \over n_0} \sZ} q^n R_n \in K_{S^1}(M)[[q^{1 \over n_0}]]$,
 then $P$ acts on $R(q)$ by multiplication by $n$ on $R_n$. 
{}From now on, we simply write ${\rm Sym} _{q^n}(TX),\ \Lambda_{q^n}(V)$ 
as  ${\rm Sym} (TX_n),\ \Lambda (V_n)$ respectively. In this way, 
$P$ acts on $TX_n$ and $V_n$ by multiplication by $n$, and the actions $P$ on 
${\rm Sym} (TX_n),\ \Lambda (V_n)$  are  naturally induced from its actions on $TX_n$ and $V_n$, etc. 
So the eigenspace of $P=n$ is given by the coefficient of $q^n$ 
of the corresponding element $R(q)$. 
For $R(q)= \oplus_{n\in {1 \over n_0} \sZ} q^n R_n 
\in K_{S^1}(M)[[q^{1 \over n_0}]]$, we will also write 
\begin{eqnarray}
{\rm Ind} (D^X \otimes R(q), m , h) = {\rm Ind} (D^X \otimes R_m, h).
\end{eqnarray}

Let $H$ be the canonical basis of ${\rm Lie} (S^1) = \bR$, 
i.e., $ \exp (t H) = {\exp}  (2 \pi i t)$ for $t\in \bR$.  
If $E$ is an $S^1$-equivariant vector bundle over $M$,
let $L_H$ denote the correponding infinitesimal action of $H$ on $\Gamma (M, E)$.
 
On the fixed point set $F$, let $J_H$ be the representation of 
${\rm Lie} (S^1)$ on $E|_F$. Then on $\Gamma (F, E|_F)$, $L_H$ is  
exactly the operator $J_H$ on $\Gamma(F, E|_F)$, and the weight of the $S^1$ action 
on $\Gamma(F,E|_F)$ is given by the action 
\begin{eqnarray}
\bJ_H= {-1 \over 2 \pi} \sqrt{-1} J_H.
\end{eqnarray}

Recall that the $\bZ_2$ grading on 
$S(TX) \otimes _{n=1}^\infty {\rm Sym} (TX_n)$ 
(resp. $S(TY, \otimes_{0<v} (\det N_v)^{-1}) \otimes {\cal F}^{-p} (X)$) 
is induced by the $\bZ_2$-grading on $S(TX)$
(resp. $S(TY, \otimes_{0<v} (\det N_v)^{-1}) $). Write
\begin{eqnarray}\begin{array}{l}
F^1_V = S(V)  \bigotimes _{n=1}^\infty \Lambda (V_n) ,\\
F^2_V= \otimes_{n\in \bN + {1 \over 2}} \Lambda (V_n) .
\end{array}\end{eqnarray}
There are two natural $\bZ_2$ gradings on $F^1_V,\ F^2_V$. 
The first grading is induced by the
$\bZ_2$-grading of $S(V)$ and the nature ${\bf Z}_2$ grading 
(induced by forms of homogeneous degree ) of 
$\bigotimes _{n=1}^\infty \Lambda (V_n)$ and  
$\otimes_{n\in \bN + {1 \over 2}} \Lambda (V_n) $. 
We define $\tau_{e|F_V^{i\pm}}= \pm 1$ to be  the involution defined by this 
$\bZ_2$-grading. 
The second grading is the one
for which $F^i_V$ $ (i=1,\ 2)$ are purely even, 
i.e., $F^{i+}_V=F^i_V$. We denote by $\tau_s= {\rm Id}$ the involution 
defined by  this $\bZ_2$ grading. Then the coefficient of 
$q^n$ $ (n\in {1 \over 2} \bZ)$ in (2.13) of $R_1(V),\ R_2(V)$ 
(resp. $R_3(V),\ R_4(V)$) is exactly the  $\bZ_2$-graded vector sub-bundle of 
$(F^1_V, \tau_s)$, $(F^1_V, \tau_e)$ 
(resp. $ (F^2_V, \tau_e)$, $(F^2_V, \tau_s)$), 
on which $P$ acts  by multiplication by $n$.

We will denote by $\tau_e$ (resp. $\tau_s$) the $\bZ_2$-grading on 
$S(TX) \otimes _{n=1}^\infty {\rm Sym} (TX_n) \otimes F^i_V$ 
induced by the above $\bZ_2$-gradings.

Let $h^{V_v}$ be the metric on $V_v$ induced by the  metric $h^V$ on $V$.
In the following, we will identify $\Lambda V_v$ with 
$\Lambda \overline{V}^*_v$ by using the Hermitian metric $h^{V_v}$.
By (\ref{hyp3}), 
as in (\ref{spinor1}), there is a natural isomorphism between 
$\bZ_2$-graded $C(V)$-Clifford modules over $F$,
\begin{eqnarray}\label{f6}
S(V_0^{\bf R}, \otimes_{0<v} (\det V_v)^{-1}) \otimes_{0<v} \Lambda 
V_v \simeq S(V)_{|F}.
\end{eqnarray}

By using above notations, on the fixed point set $F$, we rewrite (\ref{hyp8}), 
for $p \in \bN$,
\begin{eqnarray}\label{f1}\begin{array}{l}
{\cal F}_p (X) = \bigotimes_{0<v} \Big (\bigotimes_{n=1}^\infty 
{\rm Sym} (N_{v,n}) 
\bigotimes_{\stackrel{n\in \bN,}{ n>pv}} {\rm Sym} (\overline{N}_{v,n}) \Big ) 
\bigotimes_{n=1}^\infty {\rm Sym} (TY_n),\\
 {\cal F}'_p (X) = \bigotimes_{ \stackrel{0<v, n\in \bN,}{0\leq n \leq pv}}
\Big (  {\rm Sym} (N_{v,-n})  \otimes \det N_v \Big ),\\
{\cal F}^{-p}(X) = {\cal F}_p (X) \otimes {\cal F}'_p (X).
\end{array}\end{eqnarray}

 Let $V_0= V_0^\sR \otimes _\sR \bC$. From  (\ref{hyp2}),  (\ref{f6}), we get
\begin{eqnarray}\label{f2}\begin{array}{l}
{\cal F}^0(X) = \bigotimes_{n=1}^\infty {\rm Sym} \Big (\oplus_{0<v}( N_{v,n} 
\oplus \overline{N}_{v,n} ) \Big ) \bigotimes_{n=1}^\infty {\rm Sym} (TY_n) \\
\hspace*{20mm} \bigotimes  {\rm Sym} (\oplus_{0<v} N_{v,0}) \otimes 
\det(\oplus_{0<v} N_{v}),\\
F^1_V =   \bigotimes_{n=1}^\infty \Lambda (\oplus_{0<v}( V_{v,n} \oplus 
\overline{V}_{v,n} ) \oplus V_{0,n}) \\
\hspace*{10mm} \otimes S(V_0^{\sR}, \otimes_{0<v} (\det V_v)^{-1})
 \otimes _{0<v}\Lambda ( V_{v,0}),\\
F^2_V =   \bigotimes_{0<n\in \sZ + 1/2} \Lambda (\oplus_{0<v}( V_{v,n} \oplus 
\overline{V}_{v,n}) \oplus V_{0,n}) . 
\end{array}\end{eqnarray}

Now,  we can reformulate Theorem 2.6 as the following Theorem.
\begin{thm}
 For each $\alpha$,  $h,\ p\in \bZ$, $p>0$, 
$m\in {1\over 2} \bZ$,  for $i=1,\ 2$, $\tau = \tau_e$ or $\tau_s$, 
 we have the following identity in $K(B)$,
\begin{eqnarray}\label{f3}\begin{array}{l}
{\rm Ind}_\tau (D^{Y_\alpha} \otimes {\cal F}^{-p} (X)  \otimes F^i_V ,
m + {1\over 2} p^2 e(N) +{1\over 2} p d'(N), h)\\
={\rm Ind}_\tau (D^{Y_\alpha} \otimes {\cal F}^0 (X)  \otimes F^i_V,
m + ph + p^2 e, h+ 2p e ).
\end{array}\end{eqnarray}
\end{thm}

{\em Proof} : The rest of this section  is devoted to the proof of  Theorem 3.1. 
\hfill $\blacksquare$

\subsection{ \normalsize Proof of Theorem 3.1}

Inspired by [{\bf T}, \S 7], for $p\in \bN^*$, we define the shift operators, 
\begin{eqnarray}\begin{array}{l}
r_*: N_{v, n} \to N_{v, n+pv}, \qquad 
r_*: \overline{N}_{v, n} \to \overline{N}_{v, n-pv}, \\
r_*: V_{v, n} \to V_{v, n+pv}, \qquad 
r_*: \overline{V}_{v, n} \to \overline{V}_{v, n-pv}. 
\end{array}\end{eqnarray}
This means that we change the action of the operator $P$ on 
$N_{v,n}$ (resp. $\overline{N}_{v,n}$) by $n+pv$ (resp. $n-pv$), etc.

Recall that $L(N),\ L(V)$ are the complex line bundles over $F$ 
defined by (\ref{hyp5}).
Also recall that $L= L(N)^{-1} \otimes L(V)$ is a trivial complex line bundle 
on $F$, and $g\in S^1$ acts on it by multiplication by $g^{2e}$.

\begin{prop} For $p\in \bZ$, $p>0$, $i=1,\ 2$, 
 there are  natural isomorphisms of vector bundles over $F$,
\begin{eqnarray}\label{shift1}
\begin{array}{l}
r_* ({\cal F}^{-p} (X)) \simeq {\cal F}^{0} (X) \otimes L(N)^p,\\
r_* (F^i_V) \simeq F^i_V \otimes L(V)^{-p}.
\end{array}\end{eqnarray}
\end{prop}

{\em Proof} : 
1) Under the action of the shift operator $r_*$,
\begin{eqnarray}\label{FX}
\begin{array}{l}
r_* ({\cal F}'_{p} (X)) = \bigotimes_{ \stackrel{0<v}{0\leq n \leq pv}}
\Big (  {\rm Sym} (N_{v,-n+pv})  \otimes \det N_v \Big )\\
=  \bigotimes_{ \stackrel{0<v}{0\leq n \leq p v}}
  {\rm Sym} (N_{v,n})\otimes_{0<v} \det N_v \otimes L(N)^p,
\end{array}\end{eqnarray}
{}From (3.5), (\ref{FX}), we get the first equation of (\ref{shift1}).

2) For $F^i_V$ $(i=1,\ 2)$, we only need to consider the shift operator on the following 
elements
\begin{eqnarray}\label{FV1}
\begin{array}{l}
F^1_{V,F} =   \bigotimes_{n=1}^\infty \Lambda (\oplus_{0<v}( V_{v,n} \oplus 
\overline{V}_{v,n} ) ) 
\otimes _{0<v}\Lambda ( V_{v,0}) ,\\
F^2_{V,F} =   \bigotimes_{n\in \bN + 1/2} \Lambda (\oplus_{0<v}( V_{v,n} \oplus 
\overline{V}_{v,n}) ) .
\end{array}\end{eqnarray}
We compute easily that
\begin{eqnarray}\label{FV2}
\begin{array}{l}
r_* F^1_{V,F}  =   \bigotimes_{n=1}^\infty 
\Lambda (\oplus_{0<v}( V_{v,n+pv} \oplus \overline{V}_{v,n-pv} ) ) 
\otimes _{0<v}\Lambda ( V_{v,pv}) ,\\
r_* F^2_{V,F}  =   \bigotimes_{n\in \bN + 1/2} 
\Lambda (\oplus_{0<v}( V_{v,n+pv} \oplus \overline{V}_{v,n-pv}) ) .
\end{array}\end{eqnarray}
The Hermitian metric $h^{V_v}$ on $V_v$ induces a natural isomorphism 
of complex  vector bundles over $F$,
\begin{eqnarray}\label{isomV}
\Lambda ^i \overline{V}_{v} \simeq \Lambda ^{\dim V_v -i} V_{v} 
\otimes \det  \overline{V}_v. 
\end{eqnarray}
In fact, let $dv_{V_v}$ be the volume form on $(V_{v,\sR}, h^{V_v})$, then 
we define 
$\Phi: \Lambda ^{\dim V_v -i} V_{v} 
\otimes \det  \overline{V}_v \to (\Lambda ^i V_v)^*$ as follows: 
for $s_1 \in \Lambda ^{\dim V_v -i} V_{v} 
\otimes \det  \overline{V}_v,\ s_2 \in \Lambda ^i V_v$,
\begin{eqnarray}
\Phi(s_1) (s_2) dv_{V_v}= s_1 \wedge s_2. \nonumber 
\end{eqnarray}
Clearly, $\Phi$ is an isomorphism of complex vector bundles. 
By using the Hermitian metric $h^{V_v}$, we identify $(\Lambda ^i V_v)^*$ to 
$\Lambda ^i \overline{V}_v$.

For $n\in \bN,\ 0<  n \leq pv$, $0\leq i \leq \dim V_v$, 
(\ref{isomV}) induces a natural $S^1$-equivariant 
 isomorphism of complex vector bundles
\begin{eqnarray}\label{FV3}
\begin{array}{l}
\Lambda ^i \overline{V}_{v,n-pv} \simeq \Lambda ^{\dim V_v -i} V_{v, -n+pv} 
\otimes \det  \overline{V}_v,\\
\Lambda ^i \overline{V}_{v,n-pv- {1 \over 2}}
\simeq \Lambda ^{\dim V_v -i} V_{v, -n+pv+ {1 \over 2}} 
\otimes \det  \overline{V}_v. 
\end{array}\end{eqnarray}
This means 
\begin{eqnarray}\label{FV4}
\begin{array}{l}
\displaystyle{
\bigotimes_{n\in \bN, 0< n \leq p v} \Lambda ^{i_n} \overline{V}_{v,n-pv} \simeq 
\bigotimes_{n\in \bN,0< n \leq pv} \Big  (\Lambda ^{\dim V_v -i_n} V_{v, -n+pv} 
\otimes \det  \overline{V}_v \Big ),    }\\ 
\displaystyle{
\bigotimes_{n\in \bN,0\leq n < pv} 
\Lambda ^{i_n} \overline{V}_{v,n-pv+ {1 \over 2}} 
\simeq \bigotimes_{n\in \bN,0\leq n < pv} 
\Big  (\Lambda ^{\dim V_v -i_n} V_{v, -n+pv- {1 \over 2}} 
\otimes \det  \overline{V}_v \Big ).  } 
\end{array}\end{eqnarray}
{}From (2.10) and the isomorphisms (\ref{FV2}) and (\ref{FV4}) 
of complex vector bundles over $F$, one gets the second induced isomorphism in
(\ref{shift1}).

The proof of Proposition 3.1 is complete. \hfill $\blacksquare$

\begin{prop} For $ p\in \bZ$, $p>0$,  $i=1,\ 2$, 
the bundle isomorphism induced by (\ref{shift1}),
\begin{eqnarray}\label{tran0}\begin{array}{l}
r_* : S(TY, \otimes_{0<v} (\det N_v)^{-1}) \otimes {\cal F}^{-p} (X) 
\otimes F^i_V    \\
\hspace*{10mm}\to S(TY, \otimes_{0<v} (\det N_v)^{-1})
\otimes {\cal F}^{0} (X) \otimes
F^i_V  \otimes L^{-p},
\end{array}\end{eqnarray}
verifies the following identities:
\begin{eqnarray}\label{tran1}
\begin{array}{l}
r_*^{-1}\cdot \bJ_H\cdot r_* =   \bJ_H,\\
r_*^{-1}\cdot P\cdot r_* =P +  p\bJ_H 
+p^2 e - {1 \over 2} p^2  e(N) - {p \over 2} d'(N).
\end{array}\end{eqnarray}
For the $\bZ_2$-gradings, we have
\begin{eqnarray}\label{tran2}
\begin{array}{l}
\tau_e r_* = r_* \tau_e,\qquad
\tau_s r_* = r_* \tau_s.
\end{array}\end{eqnarray}
\end{prop}

{\em Proof} : The first 
equation of (\ref{tran1}) is obvious.

To prove the second equation of (\ref{tran1}) 
we divide the argument into several steps.

 a) From (\ref{FV4}),
on $\bigotimes_{\stackrel{n\in \bN}{0< n \leq pv}} 
\Lambda ^{i_n} \overline{V}_{v,n}$, we have 
\begin{eqnarray}\label{tran3}\begin{array}{l}
r_*^{-1} P r_* = \sum _{\stackrel{n\in \bN}{0< n \leq pv}} 
(\dim V_v - i_n) (-n + pv) \\
\hspace*{20mm}
= P + p\bJ_H +\sum _{\stackrel{n\in \bN}{0< n \leq pv}}(-n + pv)\dim V_v \\
\hspace*{20mm} = P + p\bJ_H + {1 \over 2} (p^2v^2 - p v) \dim V_v.
\end{array}\end{eqnarray}
Thus,  from (\ref{hyp4}), (\ref{tran3}),  on 
$\bigotimes_{0<v} \bigotimes_{\stackrel{n\in \bN}{0< n \leq pv}} 
\Lambda ^{i_n} \overline{V}_{v,n}$,
we have
\begin{eqnarray}\label{tranP1}
r_*^{-1} P r_* = P+ p\bJ_H + {1 \over 2} p^2 e(V) - {1 \over 2} pd'(V).
\end{eqnarray}

The operators $P$ and $\bJ_H$ act on $S \Big (V_0^{\sR}, \det (\oplus_{0<v} 
V_v)^{-1} \Big )$ 
by multiplication by 
\begin{eqnarray}\label{tranP2}
0,\quad   
-{1 \over 2} \Sigma_{0<v} v \dim V_v= -{1\over 2} d'(V).
\end{eqnarray}
respctively.

When $P$ acts the rest part of $F^1_V$,
we have 
\begin{eqnarray}\label{tranP3}
r_*^{-1} P r_* = P+ p\bJ_H. 
\end{eqnarray}
{}From (\ref{f2}), (\ref{tranP1}), (\ref{tranP2}) and  (\ref{tranP3}), we know that when acting on 
 $F^1_V$, one has the equality:
\begin{eqnarray}\label{tranP4}
r_*^{-1} P r_* = P+ p\bJ_H+ {1 \over 2}p^2 e(V). 
\end{eqnarray}

b) Similar to (\ref{tranP1}), from (\ref{FV4}),  on
 $\bigotimes_{\stackrel{n\in \bN}{0\leq n < pv}} \Lambda ^{i_n} 
\overline{V}_{v,n+ {1 \over 2}}$, we have 
 \begin{eqnarray}\label{tranP5}
\qquad \begin{array}{l}
r_*^{-1} P r_* = \sum _{\stackrel{n\in \bN}{0\leq n <p v}} 
(\dim V_v - i_n) (-n + pv- {1 \over 2}) \\
= P + p\bJ_H 
+  (\dim V_v)\sum _{\stackrel{n\in \bN}{0\leq n < pv}} (-n + pv- {1 \over 2})
= P + p\bJ_H + {1 \over 2} p^2 v^2 \dim V_v.
\end{array}\end{eqnarray}

{}From (\ref{hyp4}),  (\ref{tranP5}), on 
$ \bigotimes_{0<v} \bigotimes_{\stackrel{n\in \bN}{0\leq n < pv}} \Lambda ^{i_n} 
\overline{V}_{v,n+ {1 \over 2}}$, we have 
\begin{eqnarray}\label{tranP6}
r_*^{-1} P r_* =  P + p\bJ_H + {1 \over 2}p^2 e(V).
\end{eqnarray}
When $P$ acts on the rest part of $F^2_V$, one has 
\begin{eqnarray}\label{tranP7}
r_*^{-1} P r_* = P+p \bJ_H. 
\end{eqnarray}
{}From  (\ref{tranP6}), (\ref{tranP7}), on $F^2_V$, we again have (\ref{tranP4}).

c) Note that on 
$\otimes_{0<v, 0\leq n \leq pv} \det N_v$, 
$\bJ_H$ acts as $p e(N) + d'(N)$.  By (\ref{f1}), we know that  
on ${\cal F}^{-p} (X)$,
\begin{eqnarray}
r_*^{-1} P r_* =  P + p \bJ_H -p(p e(N) + d'(N)).
\end{eqnarray}

On $S(TY, \det (\oplus_{0<v} N_v)^{-1})$, $\bJ_H$ acts as 
$- {1 \over 2} d'(N)$.
So on $S(TY, \det (\oplus_{0<v} N_v)^{-1}) 
\otimes {\cal F}^{-p} (X)$
\begin{eqnarray}\label{tranP8}
r_*^{-1} P r_* =  P + p\bJ_H - p^2 e(N) -{1 \over 2}p d'(N) .
\end{eqnarray}

{}From (2.12), (\ref{tranP4}), b) and (\ref{tranP8}), 
 we get the second equation of (\ref{tran1}).

Finally,  under our operations, the $\bZ_2$-grading $\tau_s$ does not change. 
For the $\bZ_2$-grading $\tau_e$, it changes only on  
$\bigotimes_{0<v} \bigotimes_{\stackrel{n\in \bN}{0< n \leq pv}} 
 \Lambda ^{i_n} \overline{V}_{v,n}$ of $F^1_V$ 
(resp. on $\bigotimes_{0<v} \bigotimes_{\stackrel{n\in \bN}{0\leq n < pv}}  
\Lambda ^{i_n} \overline{V}_{v,n+{1\over 2}}$
of $F^2_V$). From  (\ref{FV4}), we know that
\begin{eqnarray}\label{tranP9}
r_*^{-1} \tau_e r_*= (-1)^{\Sigma_{0<v}p v \dim V_v}  \tau_e. 
\end{eqnarray}
 As ${1\over 2}p_1(V)_{S^1}\in H^*_{S^1}(M,\bZ)$ is well-defined, from
(2.5), (2.11), we get
\begin{eqnarray}\label{tranP10}
\sum_{0<v} v \dim V_v =0 \ \sm (2),
\end{eqnarray}
{}From  (\ref{tranP9}), (\ref{tranP10}), we get (\ref{tran2}).

The proof of Proposition 3.2 is complete.
\hfill $\blacksquare$\\

{\em Proof of Theorem 3.1} : From  (\ref{hyp7}) and  Propositions  3.2, for each $\alpha$,
  $h,\ p\in \bZ$, $p>0$, 
$m\in {1\over 2} \bZ$,  and for $i=1,\ 2$, $\tau = \tau_e$ or $\tau_s$, 
we have the following indentity in $K(B)$,
 \begin{eqnarray}\label{tran11}\begin{array}{l}
 {\rm Ind}_\tau (D^{Y_\alpha} \otimes {\cal F}^{-p} (X)  \otimes F^i_V ,
 m + {1\over 2} p^2 e(N) +{1\over 2} p d'(N), h)\\
 ={\rm Ind}_\tau (D^{Y_\alpha} \otimes {\cal F}^0 (X)  \otimes F^i_V
 \otimes L^{-p} ,
 m + ph + p^2 e, h )\\
 ={\rm Ind}_\tau (D^{Y_\alpha} \otimes {\cal F}^0 (X)  \otimes F^i_V ,
 m + ph + p^2 e, h+2 pe ) .
 \end{array}\end{eqnarray}

The proof of  Theorem 3.1 is complete.
\hfill $\blacksquare$\\

\newpage

\section{ \normalsize Proof of Theorem 2.7}
\setcounter{equation}{0}

In this section, we  prove Theorem 2.7. Many arguments in this section 
are inspired by [{\bf T}, \S 6,\ 9]. We will construct a family twisted 
Dirac operator on $M(n_j)$, the fixed point set of the
induced $\bZ_{n_j}$ action on $M$. 
 By applying our $K$-theory version of the equivariant family
index theorem to this operator,
 we derive Theorem 2.7.

This section is organized as follows: 
In Section 4.1, we construct a family Dirac operator on $M(n_j)$. 
In Section 4.2, by introducing a shift operator, we will relate both sides
 of equation (\ref{e8}) to the index bundle
  of the family Dirac operator 
on $M(n_j)$. In Section 4.3, we prove Theorem 2.7.

In this section, we  make the same assumptions and use the same 
notations as in  Sections 2 and 3.

\subsection{ \normalsize The Spin$^c$ Dirac operator on $M(n_j)$}

Let $\pi: M\to B$ be a  fibration of compact manifolds with fiber $X$ 
and $\dim_{\sR} X= 2l$. We assume that  $S^1$ acts fiberwise on $M$,
and $TX$ has an $S^1$-invariant spin structure.
 Let $V$ be a real vector bundle  over $M$ carrying  an $S^1$-invariant spin 
structure and $\dim_{\sR} V= 2k$.

Let $F=\{F_\alpha \}$ be the fixed point set of the $S^1$ action on $M$. 
Then $\pi: F\to B$ is a fibration with compact fiber $Y$. For $n\in \bN,\ n>0$, let $\bZ_n \subset S^1$ denote
the cyclic subgroup of order $n$. 

For $ n_j\in \bN$ with $n_j>0$,  
let $M(n_j)$ be the fixed point set of the
induced $\bZ_{n_j}$-action on $M$. 
Then $\pi: M(n_j)\to B$ is a fibration with compact fiber $X(n_j)$.
Let $N(n_j)\to M(n_j)$ be the  normal bundle of $M(n_j)$ in $M$. 
Then we have the following 
$\bZ_{n_j}$-equivariant decomposition of $N(n_j)\otimes_\sR \bC$ over $M(n_j)$,
\begin{eqnarray}
N(n_j) \otimes_\sR \bC = \bigoplus_{0<v <n_j} N(n_j)_v.
\end{eqnarray}
Here $N(n_j)_v$ is the  complex vector bundles over $M(n_j)$ with $g\in \bZ_{n_j}$
acting by $g^v$ on it. Complex conjugation provides a $\bC$ anti-linear 
isomorphism between $N(n_j)_v$ and $\overline{N(n_j)_{n_j-v}}$. 
If $n_j$ is even, this produces a real structure on $N(n_j)_{{n_j \over 2}}$, 
so this bundle is the complexification of a real vector bundle
 $N(n_j)_{{n_j \over 2}}^{\sR}$ on $M(n_j)$.
Thus, $N(n_j)$ is isomorphic, as  a  real vector bundle, to 
\begin{eqnarray} \label{dec1}
N(n_j) \simeq  \bigoplus_{0<v <n_j/2} N(n_j)_v \oplus N(n_j)_{{n_j \over 2}}^{\sR} .
\end{eqnarray}

Similarly, we have the following 
$\bZ_{n_j}$-equivariant decomposition of $V\otimes_\sR \bC$,
\begin{eqnarray}
V \otimes_\sR \bC = \bigoplus_{0 \leq v <n_j} V(n_j)_v.
\end{eqnarray}
Here $V(n_j)_v$ is the   complex vector bundle over $M(n_j)$ with $g\in \bZ_{n_j}$
acting by $g^v$ on it. For $v\neq 0$, complex conjugation provides a 
$\bC$ anti-linear isomorphism between $V(n_j)_v$ and 
$\overline{V(n_j)_{n_j-v}}$. If $n_j$ is even, this produces a real 
structure on $V(n_j)_{{n_j \over 2}}$, so this bundle is the complexification 
of a real vector bundle $V(n_j)_{{n_j \over 2}}^{\sR}$ over $M(n_j)$. 
Complex conjugation also provides  a real structure on $V(n_j)_0$ such that 
$V(n_j)_0= V(n_j)_{0}^{\sR} \otimes_\sR \bC$. 
Thus, over $M(n_j)$,  $V$ is isomorphic, as  a  real vector bundle, to 
\begin{eqnarray}\label{dec2}
V|_{M(n_j)} \simeq V(n_j)_{0}^{\sR} \bigoplus_{0<v <n_j/2} V(n_j)_v 
\oplus V(n_j)_{{n_j \over 2}}^{\sR} .
\end{eqnarray}
In (\ref{dec1}), (\ref{dec2}), the last term is understood to be zero 
when $n_j$ is odd.

It is essential for us to know that
 the vector bundles $TX(n_j)$ and  $V(n_j)_{0}^{\sR}$ are orientable. 
For this we have
the following Lemma which was proved in [{\bf E}] and [{\bf BT}, Lemma 10.1]. 

\begin{lemma}  Let $W$ be a real, spin vector bundle over a manifold $M$.
 We assume that $\bZ_n$ $(n\in \bN^*)$ acts on $M$,
 and that the $\bZ_n $ action lifts on $W$ and preserves the spin structure of $W$.
 Let  $M(n)$ be the fixed point set of the $\bZ_n$ action on $M$. Let $W_0$ be the 
subbundle of $W$ over $M(n)$ on which the generator of $\bZ_n$ acts trivially.
Then $W_0$ is orientable.
\end{lemma}

By  Lemma 4.1,  $TX(n_j)$ and $V(n_j)_{0}^{\sR}$ are orientable over $M(n_j)$. 
Thus $N(n_j)$ is   orientable over $M(n_j)$. By (\ref{dec1}), (\ref{dec2}),
$N(n_j)_{{n_j \over 2}}^{\sR}$ and $V(n_j)_{{n_j \over 2}}^{\sR}$ 
are also orientable over $M(n_j)$. 
In the following, we fix  the orientations of 
 $N(n_j)_{{n_j \over 2}}^{\sR}$, $V(n_j)_{{n_j \over 2}}^{\sR}$ over $M(n_j)$.
Then $TX(n_j)$ and $V(n_j)_{0}^{\sR}$  are naturally oriented by 
(\ref{dec1}), (\ref{dec2}) and the orientations of $TX$, 
$V,N(n_j)_{{n_j \over 2}}^{\sR}$, $V(n_j)_{{n_j \over 2}}^{\sR}$. 

Let us denote by 
\begin{eqnarray} \label{rn1}
r(n_j)= {1\over 2} ( 1 + (-1)^{n_j}).
\end{eqnarray}

\begin{lemma} Assume that (\ref{hyp1})  holds.  Let 
\begin{eqnarray} \label{line1}
L(n_j)= \bigotimes _{0<v <n_j/2} \Big ( \det (N(n_j)_v) \otimes 
\det (\overline {V(n_j)_v})\Big ) ^{(r(n_j)+1)v}
\end{eqnarray}
be the line bundle over $M(n_j)$. Then:

i)  $L(n_j)$ has an $n_j^{\rm th}$ root over $M(n_j)$. 

ii) Let 
  \begin{eqnarray}\label{line2}
\begin{array}{l}
L_1= \bigotimes _{0<v <n_j/2} \Big ( \det (N(n_j)_v) \otimes 
\det (\overline {V(n_j)_v})\Big ) \otimes L(n_j)^{r(n_j)/n_j},\\
L_2= \bigotimes _{0<v <n_j/2} \Big ( \det (N(n_j)_v) \Big ) 
\otimes L(n_j)^{r(n_j)/n_j}.
\end{array}\end{eqnarray}
Let $U_1= T X(n_j)\oplus V(n_j)_{0}^{\sR}$ and 
$U_2 = TX(n_j) \oplus V(n_j)_{{n_j \over 2}}^{\sR}$.
Then  $TX(n_j),\ V(n_j)_{{n_j \over 2}}^{\sR} $ and $V(n_j)_{0}^{\sR}$ 
are of even dimensions. Furthermore, $U_1$ (resp. $U_2$) has a $Spin^c$ 
structure defined by  $L_1$ (resp. $L_2$).
\end{lemma}

{\em Proof} : By [{\bf BT}, Lemma 9.4], $TX(n_j),\ V(n_j)_{{n_j \over 2}}^{\sR}$
 and $V(n_j)_{0}^{\sR}$ are of even dimensions. 
{}From the proof of [{\bf BT}, Lemmas 11.3 and 11.4],
we get the rest part of Lemma 4.2. \hfill $\blacksquare$\\

Lemma 4.2 is very important. It allows us, as we are going to see, to apply
the constructions and results in Section 1.3 to the fibration
$M(n_j)\rightarrow B$, which is the main concern of this Section.

For $p_j \in \bN$, $p_j< n_j$, $(p_j, n_j)=1$, $\beta_j = {p_j \over n_j}$, let us write 
\begin{eqnarray}\begin{array}{l}
{\cal F} (\beta_j) =  \otimes _{0<n\in \bZ} {\rm Sym} (TX(n_j)_n) 
\bigotimes_{0<v <n_j/2} {\rm Sym} 
\Big ( \bigoplus_{0< n \in \bZ + p_j v/n_j} N(n_j)_{v,n}\\
\hspace*{20mm}\bigoplus_{0< n \in \bZ - p_j v/n_j} \overline{N(n_j)}_{v,n} \Big ) 
\otimes _{0<n\in \bZ + {1 \over 2}} {\rm Sym} (N(n_j)_{{n_j \over 2},n} ),
\end{array}\label{idFV1}\\
\begin{array}{l}
F^1_V(\beta_j)  = \Lambda \Big (\oplus_{0<n\in \bZ} V(n_j)_{0,n} \bigoplus_{0< v <n_j/2} 
\Big ( \oplus_{0<n \in \bZ + p_j v/n_j} V(n_j)_{v,n} \\
\hspace*{20mm}\oplus _{0< n \in \bZ - p_j v/n_j} \overline{V(n_j)}_{v,n} \Big ) 
\oplus _{0< n\in \bZ + {1 \over 2} } V(n_j)_{{n_j \over 2}, n} \Big ),\\
F^2_V(\beta_j) = \Lambda \Big (\oplus_{0<n\in \bZ} V(n_j)_{{n_j \over 2},n} 
\bigoplus_{0< v <n_j/2} \Big ( \oplus_{0<n \in \bZ + p_j v/n_j+ {1 \over 2}} 
V(n_j)_{v,n}    \\
\hspace*{20mm}\oplus _{0< n \in \bZ - p_j v/n_j +{1 \over 2}} 
\overline{V(n_j)}_{v,n}\Big ) \oplus _{0< n\in \bZ + {1 \over 2} } 
V(n_j)_{0, n} \Big ).
\end{array}  \nonumber 
\end{eqnarray}
 
We denote by $D^{X(n_j)}$ the $S^1$-equivariant Spin$^c$-Dirac operator 
on $S(U_1,L_1)$  or $S(U_2, L_2)$ along the fiber $X(n_j)$ 
defined as in Section 1.3. 
We denote by $D^{X(n_j)} \otimes {\cal F}(\beta_j)
\otimes F^i_V(\beta_j)$ $(i=1,\ 2)$ the corresponding twisted Spin$^c$ Dirac 
operator on $S(U_i,L_i) \otimes {\cal F}(\beta_j)
\otimes F^i_V(\beta_j)$ along the fiber $X(n_j)$. \\

{\bf Remark 4.1.}  In fact, to define an $S^1$-action on $ L(n_j)^{r(n_j)/n_j}$, 
one must replace the $S^1$-action by its $n_j$-fold action. Here by abusing 
notation, we still say an $S^1$-action without causing any confusion.\\

In the rest of this subsection, we will reinterpret all of the above objects
when we restrict ourselves to $F$, the fixed point set of the $S^1$ action.
We will use the notation of Sections 1.3,  2.

Let $N_{F/M(n_j)}$ be the normal bundle of $F$ in $M(n_j)$. Then by (2.7),
\begin{eqnarray}\label{idFV2}\begin{array}{l}
N_{F/M(n_j)}= \bigoplus _{0<v : v\in n_j \bZ} N_v,\\
TX(n_j) \otimes_\sR  \bC = TY\otimes_\sR \bC \oplus_{0<v, v\in n_j \bZ}
 (N_v \oplus \overline{N}_v).
\end{array}\end{eqnarray}
By (\ref{hyp2}) and (4.1), the restriction of 
$N(n_j)_v $ $(1\leq v \leq n_j/2)$  to $F$ is given by 
\begin{eqnarray}\label{idFV3}
N(n_j)_v = \bigoplus _{0<v': v'=v\ \sm (n_j) } N_{v'} 
\bigoplus _{0<v': v'= -v\ \sm  (n_j) } \overline{N}_{v'}.
\end{eqnarray}
By (\ref{hyp3}) and (\ref{dec2}), the restriction  of 
$V(n_j)_v$ $ (1\leq v \leq n_j/2)$ to $F$ is given by 
\begin{eqnarray}\label{idFV4}
V(n_j)_v = \bigoplus _{0<v': v'=v\ \sm (n_j) } V_{v'} 
\bigoplus _{0<v': v'= -v\ \sm  (n_j) } \overline{V}_{v'}.
\end{eqnarray}
and for $v=0$,
\begin{eqnarray}\label{idFV5}
V(n_j)_0 = V_{0}^{\sR}  \otimes_\sR \bC  \bigoplus_{0<v, v=0\ \sm (n_j)}
 (V_v \oplus \overline{V}_v).
\end{eqnarray}
{}From  (\ref{idFV2})-(\ref{idFV5}), 
we have the following identifications of  real vector bundles over $F$,
\begin{eqnarray}\label{idFV6}\begin{array}{l}
N(n_j)_{{n_j \over 2}}^{\sR} =  \bigoplus_{0<v, v={n_j \over 2}\ \sm (n_j)}
 {N}_v ,\\
TX(n_j) = TY \bigoplus_{0<v, v=0\ \sm (n_j)}
 {N}_v ,\\
V(n_j)_{0}^{\sR} = V_{0}^{\sR} \bigoplus_{0<v, v=0\ \sm (n_j)}
 V_v ,\\
V(n_j)_{{n_j \over 2}}^{\sR} = \bigoplus_{0<v, v={n_j \over 2}\ \sm (n_j)}
  V_v.
\end{array}\end{eqnarray}

We denote by $V_0= V_{0}^{\sR} \otimes_\sR \bC$ the complexification of
 $V_{0}^{\sR}$ over $F$.  As $(p_j, n_j)=1$, we know that, for $v\in\bZ$, 
$p_jv /n_j\in \bZ$ iff $v/n_j\in \bZ$. 
Also, $p_jv /n_j\in \bZ+ {1 \over 2}$ iff $v/n_j\in  \bZ+ {1 \over 2}$.
 From (\ref{idFV2})-({\ref{idFV5}), we then get
\begin{eqnarray}\label{idFV7}\\
\begin{array}{l}
{\cal F} (\beta_j) = \otimes_{0<n\in \sZ} {\rm Sym} (TY_n)
\bigotimes_{0<v, v=0, {n_j \over 2}\ \sm (n_j) } 
\bigotimes_{0<n\in \bZ + {p_j v \over n_j}} 
{\rm Sym} (N_{v,n} \oplus \overline{N}_{v,n}) \\
\hspace*{10mm}  \bigotimes_{0<v'< n_j/2}
{\rm Sym} \Big ( \oplus_{v=v'\ \sm (n_j)} \Big ( \oplus_{0<n\in \bZ 
+ {p_j v \over n_j}} N_{v,n} 
\oplus_{0<n\in \bZ - {p_j v \over n_j}} \overline{N}_{v,n} \Big ) \\
\hspace*{15mm} \oplus_{v=-v'\ \sm (n_j)} \Big (  \oplus_{0<n\in \bZ 
+ {p_j v \over n_j}} N_{v,n} 
\oplus_{0<n\in \bZ - {p_j v \over n_j}} \overline{N}_{v,n} \Big ) \Big ),\\
\\
F^1_V(\beta_j) = \Lambda  \Big [ \oplus_{0< n\in \bZ} V_{0,n} 
\bigoplus_{0<v, v=0, {n_j \over 2}\ \sm (n_j)} 
\Big ( \oplus_{0< n\in \bZ +  {p_j v \over n_j}} V_{v,n} 
 \oplus_{0< n\in \bZ -  {p_j v \over n_j}} \overline{V}_{v,n} \Big )\\
\hspace*{10mm}  \bigoplus_{0<v' < n_j/2} 
\Big ( \bigoplus_{v=v', -v'\ \sm (n_j)} 
\Big (  \oplus_{0< n\in \bZ +  {p_j v \over n_j}} V_{v,n} 
 \oplus_{0<n\in \bZ - {p_j v \over n_j}} \overline{V}_{v,n} 
\Big ) \Big )  \Big ],\\
F^2_V(\beta_j) = \Lambda \Big [\oplus_{0< n\in \bZ+ {1 \over 2}} V_{0,n} 
\oplus_{0<v, v=0, {n_j \over 2}\ \sm (n_j) } 
\Big ( \oplus_{0< n\in \bZ +  {p_j v \over n_j}+ {1 \over 2}} V_{v,n}  
 \oplus_{0< n\in \bZ -  {p_j v \over n_j}+ {1 \over 2}} \overline{V}_{v,n} 
\Big )\\
\hspace*{10mm}  \bigoplus_{0<v' < n_j/2} \Big ( \oplus_{v=v', -v'\ \sm (n_j)} 
\Big (  \oplus_{0< n\in \bZ +  {p_j v \over n_j}+ {1 \over 2}} V_{v,n}  
\oplus_{0<n\in \bZ - {p_j v \over n_j}+ {1 \over 2}} \overline{V}_{v,n} 
\Big ) \Big ) \Big ].
\end{array} \nonumber
\end{eqnarray}

Now we want to compare the spinor bundles  over $F$. From (\ref{line1}), (\ref{line2}), (\ref{idFV3}) and  (\ref{idFV4}),  
we find  that over $F$ we have 
\begin{eqnarray} \label{line3} \begin{array}{l}
L(n_j)^{r(n_j) \over n_j} = \bigotimes_{0< v' < n_j/2} \Big ( 
\otimes _{v=v' \ \sm (n_j)} (\det N_v \otimes \det \overline{V}_v) ^{2v'}   \\
\hspace*{25mm} \otimes _{v=-v'\ \sm (n_j)} 
(\det N_v \otimes \det \overline{V}_v) ^{-2v'} \Big )^{r(n_j)/n_j},\\
L_1= \bigotimes_{0< v' < n_j/2} \Big ( 
\otimes _{v=v' \ \sm (n_j)} (\det N_v \otimes \det \overline{V}_v)    \\
\hspace*{25mm}\otimes _{v=-v'\ \sm (n_j)} 
(\det N_v \otimes \det \overline{V}_v) ^{-1} \Big )
  \otimes L(n_j)^{r(n_j)/n_j},\\
L_2= \bigotimes_{0< v' < n_j/2} \Big ( 
\otimes _{v=v'\  \sm (n_j)} \det N_v   \\
\hspace*{25mm} \otimes _{v=-v'\ \sm (n_j)} (\det N_v ) ^{-1} \Big )  
\otimes L(n_j)^{r(n_j)/n_j}.
\end{array}\end{eqnarray}

{}From  (\ref{idFV6}), over $F$, we have
\begin{eqnarray} \label{idFV8} \begin{array}{l}
TX(n_j) \oplus V(n_j)_{0}^{\sR} = TY \oplus V_{0}^{\sR} 
  \oplus _{0<v, v=0\ \sm (n_j)} (N_v \oplus V_v),\\
TX(n_j) \oplus V(n_j)_{{n_j \over 2}}^{\sR}  
= TY  \oplus _{0<v, v=0 \ \sm (n_j)} N_v 
\oplus_{0<v, v={n_j \over 2}\ \sm (n_j)}  V_v.
\end{array}\end{eqnarray}
Recall that the Spin$^c$ vector bundles $U_1$ and $U_2$ have been defined in
Lemma 4.2. Let us write 
\begin{eqnarray} \label{idFV11} \\
\begin{array}{l}
\displaystyle{S(U_1,L_1)' = S \Big (TY\oplus V_{0}^{\sR}, 
L_1 \bigotimes_{\stackrel{0<v,}{v=0\ \sm (n_j)}} 
(\det N_v \otimes \det V_v)^{-1}  \Big ) 
\bigotimes _{\stackrel{0<v,}{v=0\ \sm (n_j)}} \Lambda V_v,  }\\
\displaystyle{S(U_2,L_2)' = S \Big (TY, L_2 \bigotimes_{\stackrel{0<v,}{v=0\ \sm (n_j)}} 
(\det N_v)^{-1} \bigotimes _{\stackrel{0<v,}{v={n_j \over 2}\ \sm (n_j)}}
( \det V_v)^{-1} \Big  ) \bigotimes _{\stackrel{0<v,}{v={n_j \over 2} 
\ \sm (n_j)}} \Lambda V_v.  }
\end{array}\nonumber
\end{eqnarray}
Then from  (\ref{dirac3}), (1.49) and  
(\ref{idFV11}), for $i=1,\ 2$, we have the following isomorphism of 
Clifford modules over $F$,
\begin{eqnarray} \label{idFV12} \begin{array}{l}
S(U_i,L_i) \simeq  S(U_i,L_i)' 
\otimes \Lambda (\oplus _{0<v, v=0\ \sm (n_j)} N_v).
\end{array}\end{eqnarray}
We define the $\bZ_2$ gradings on $S(U_i,L_i)'\ (i=1,\ 2)$ as that induced by the 
$\bZ_2$-gradings on $S(U_i,L_i)$ $(i=1,\ 2)$ and on 
$\Lambda (\oplus _{0<v, v=0 \ \sm (n_j)} N_v)$ such that the isomorphism
(\ref{idFV12}) preserves these $\bZ_2$-gradings.

We define formally the following complex line bundles over $F$,
\begin{eqnarray}   \qquad\begin{array}{l}
L'_1= \Big [ L_1^{-1} \otimes _{\stackrel{0<v,}{ v=0\ \sm (n_j)}} 
(\det N_v \otimes \det V_v)\otimes _{0<v} 
(\det N_v \otimes \det V_v)^{-1}  \Big ] ^{1/2}, \\
L'_2= \Big [ L_2^{-1}  \otimes _{\stackrel{0<v,}{ v=0\ \sm (n_j)}} 
\det N_v \otimes _{\stackrel{0<v,}{ v=n_j/2\ \sm (n_j)}}  \det V_v
\otimes _{0< v} (\det N_v)^{-1}   \Big ]^{1/2}.
\end{array} \nonumber\end{eqnarray}
{}From (\ref{chern1}),  (\ref{dirac3}), Lemma 4.2 and the assumption that $V$ is
spin, one verifies easily that 
$c_1({L'_i}^2) = 0 \ \sm (2)$ for $i=1,\  2$.
Thus $L'_1,\ L'_2$ are well defined complex line bundles over $F$.
For later use, we also write down the following expressions of
$L'_i$ ($i=1,\ 2$) which can be deduced from (\ref{line3}):
\begin{eqnarray}  \label{line4} \begin{array}{l}
L'_1= \Big [ L(n_j)^{-1/n_j} \otimes _{0<v, v={n_j \over 2}\ \sm (n_j)} 
(\det N_v \otimes \det \overline{V}_v) \Big ] ^{{r(n_j) \over 2}} \\
\hspace*{20mm}\otimes _{0< v, 0< v \leq {n_j \over 2}\ \sm (n_j)} 
(\det N_v)^{-1} 
\otimes _{0< v, {n_j \over 2}< v <n_j \ \sm (n_j)}  (\det V_v )^{-1},\\
L'_2= \Big [ L(n_j)^{-1/n_j} \otimes _{0<v, v={n_j \over 2}\ \sm (n_j)} 
( \det N_v  \otimes \det V_v) 
\Big ] ^{{r (n_j)\over 2}}\\
\hspace*{20mm}
\otimes _{0< v, 0 < v \leq {n_j \over 2} \sm (n_j)} (\det N_v)^{-1}.
\end{array}\end{eqnarray}

{}From (\ref{line3}), (\ref{idFV11}) and  the definition of $L'_i$ $(i=1,\ 2)$,
we get the following identifications of Clifford modules over $F$,
\begin{eqnarray} \label{idFV9}\qquad  \begin{array}{l}
S(U_1, L_1 )' \otimes L'_1 = S(TY, \otimes _{0<v} (\det N_v)^{-1})
\otimes S(V_{0}^{\sR} , \otimes _{0<v} (\det V_v)^{-1})\\
\hspace*{30mm}\otimes \Lambda (\oplus _{0<v, v=0 \ \sm (n_j)} V_v),\\
S(U_2, L_2 )' \otimes L'_2 = S(TY, \otimes _{0<v} (\det N_v)^{-1})
\otimes \Lambda (\oplus_{0<v, v={n_j \over 2} \ \sm (n_j)} V_v).
\end{array}\end{eqnarray}

To compare the $\bZ_2$-gradings in (\ref{idFV9}), 
we will compare explicitly the orientations. Recall that, if $(W, h^W)$ is a real Euclidian vector space of dimension $2m$, and $J$ is a complex structure on $W$ which preserves $h^W$. 
Let $\{e_i, J e_i\}_{i=1}^m$ be an orthonormal basis of $(W, h^W)$. 
Then $W$ is canonically oriented, and its orientation is defined by 
the canonical Riemannian volume form 
\begin{eqnarray}
(e_1\wedge J e_1) \wedge \cdots (e_m \wedge J e_m) = d v _W. \nonumber
\end{eqnarray}

Let $dv_{TX},\  dv_V$ be the corresponding Riemannian volume forms on $(TX, h^{TX})$ and  
$(V, h^V)$ which define the orientations of $TX,\  V$  over $M$. 
Let $dv_{N_v}$ and $dv_{\overline{N}_v}$ ( resp. $dv_{V_v}$, 
$dv_{\overline{V}_v}$) $(0<v)$ be the canonical Riemannian volume forms
 on $N_{v, \sR}$ and $\overline{N}_{v, \sR}= N_{v, \sR}$ 
(resp. $V_{v, \sR}$,  $\overline{V}_{v, \sR}= V_{v, \sR}$). 
Then through the identifications (\ref{spinor1}) and (\ref{f6}), 
 the orientations of $TY$ and $V_{0}^{\sR} $  over $F$ are 
defined by the volume forms $dv_{TY}$ and $dv_{V_{0}^{\sR}}$ respectively such that 
\begin{eqnarray}\label{vol1}\begin{array}{l}
dv_{TX}= dv_{TY} \otimes _{0<v} dv_{N_v},\\
dv_{V}= dv_{V_{0}^{\sR} } \otimes _{0<v} dv_{V_v}. 
\end{array}\end{eqnarray}

By (\ref{dec1}), (\ref{dec2}), (4.10) and (4.11),
the orientations of $TX(n_j)$ and  
$V(n_j)_{0}^{\sR} $, when restricted to $F$, are given by 
\begin{eqnarray}\label{vol2}
\\
dv_{TX}|_F = dv_{TX(n_j)} \bigotimes _{0< v' < n_j/2} 
\Big ( \bigotimes_{v=v'\ \sm (n_j)} dv_{N_v} \bigotimes_{v=-v'\ \sm (n_j)}
 dv_{\overline{N}_v} \Big ) \otimes dv_{N(n_j)_{{n_j \over 2}}^{\sR} }, \nonumber\\
dv_{V}|_F = dv_{V(n_j)_{0}^{\sR} } \bigotimes _{0< v' < n_j/2} 
\Big ( \bigotimes_{v=v'\ \sm (n_j)} dv_{V_v} \bigotimes_{v=-v'\ \sm (n_j)} 
dv_{\overline{V}_v} \Big ) \otimes dv_{V(n_j)_{{n_j \over 2}}^{\sR} }. \nonumber
\end{eqnarray}
respectively.

Clearly, we have 
\begin{eqnarray}\label{vol3}\begin{array}{l}
dv_{N_v}= (-1)^{\dim N_v} dv_{\overline{N}_v},\\
dv_{V_v}= (-1)^{\dim V_v} dv_{\overline{V}_v}.
\end{array}\end{eqnarray}

{}From  (\ref{vol1}), (\ref{vol2}) and  (\ref{vol3}), we get 
\begin{eqnarray}\label{vol4}\begin{array}{l}
\displaystyle{dv_{TX(n_j)}= (-1) ^{\Delta (n_j,N)} 
dv_{TY} \bigotimes _{0<v, v =0\ \sm (n_j)} dv_{N_v},  }\\
\displaystyle{  dv_{V(n_j)_{0}^{\sR} }= 
(-1) ^{ \Delta (n_j,V)} 
dv_{V_{0}^{\sR} } \bigotimes _{0<v, v =0\ \sm (n_j)} dv_{V_v}, }
\end{array}\end{eqnarray}
where
\begin{eqnarray}\label{vol5}\begin{array}{l}
\displaystyle{
\Delta (n_j,N) = \sum_{{n_j \over 2} < v' < n_j} 
\sum _{0<v=v'\ \sm (n_j)} \dim N_v 
+ o(N(n_j)_{{n_j \over 2}}^{\sR} ),  } \\
\displaystyle{\Delta (n_j,V) = \sum_{{n_j \over 2} < v' < n_j} 
\sum _{0<v=v'\ \sm (n_j)} \dim V_v 
+ o(V(n_j)_{{n_j \over 2}}^{\sR} ), }
\end{array}\end{eqnarray}
with $o(N(n_j)_{{n_j \over 2}}^{\sR} )=  0\ {\rm or}\  1$
 (resp. $o(V(n_j)_{{n_j \over 2}}^{\sR} )= 0\ {\rm or}\  1$), 
the value depends on whether the 
given orientation on $N(n_j)_{{n_j \over 2}}^{\sR} $ 
( resp. $V(n_j)_{{n_j \over 2}}^{\sR} $) agrees or disagrees 
with the complex orientation of $\oplus_{v={n_j \over 2}\ \sm (n_j)} N_v$
(resp. $\oplus_{v={n_j \over 2}\ \sm (n_j)} V_v$).

{}From (\ref{idFV6}), (\ref{idFV12}), (\ref{vol3}) and  (\ref{vol4}),
 we see that, for the $\bZ_2$-gradings induced by $\tau_s$, the difference of the
$\bZ_2$-gradings of (\ref{idFV9}) is  $(-1)^{\Delta (n_j,N)}$; 
for the $\bZ_2$-gradings induced by $\tau_e$, the difference of the 
$\bZ_2$-gradings
of the first (resp. second ) equation of (\ref{idFV9}) is 
$(-1)^{\Delta (n_j,N)+ \Delta (n_j,V)}$ 
(resp. $(-1)^{\Delta (n_j,N)
+o(V(n_j)_{{n_j \over 2}}^{\sR} ) }$).

\subsection{\normalsize The Shift operators}

Let $p\in {\bf N}^*$ be fixed. For any $1\leq j\leq J_0$,
inspired by [{\bf T}, \S 9], 
we define the following shift operators $r_{j*}$:
\begin{eqnarray}\label{shif1}\begin{array}{l}
r_{j*}: N_{v,n} \to N_{v, n + (p-1) v + p_j v /n_j}, \quad 
r_{j*}: \overline{N}_{v,n} \to 
\overline{N}_{v, n - (p-1) v - p_j v /n_j}, \\
r_{j*}: V_{v,n} \to V_{v, n + (p-1) v + p_j v /n_j}, \quad 
r_{j*}: \overline{V}_{v,n} \to 
\overline{V}_{v, n - (p-1) v - p_j v /n_j}.
\end{array}\end{eqnarray}
This means that we change the action of the operator $P$ on $N_{v,n}$ by 
 a multiplication by $n + (p-1) v + p_j v /n_j$, etc.

If $E$  is a combination of the above bundles, we denote by
$r_{j*} E$ the bundle on which the action of $P$ is changed in 
the above way. 

{}From (\ref{e4}), (\ref{f1}),  we get
\begin{eqnarray}\label{shif2}\begin{array}{l}
{\cal F}_{p,j}(X)= {\cal F}_{p}(X) \otimes {\cal F}'_{p-1}(X) 
\bigotimes_{(v,n)\in \cup_{i=1}^{j} I^p_i} \Big({\rm Sym} ( N_{v,-n}) 
\otimes  \det N_v \Big ) \\
\hspace*{30mm}\bigotimes_{(v,n)\in \overline{I}^p_{j}} 
{\rm Sym}( \overline{N}_{v,n}) .
\end{array}\end{eqnarray}

Recal that the vector bundles $F_V^i$ $(i=1,\ 2)$ have been defined in
(3.6).

\begin{prop} There are  natural isomorphisms of vector bundles over $F$,
\begin{eqnarray}\label{shif3}
 \qquad\begin{array}{l}
r_{j*} {\cal F}_{p,j-1}(X) \simeq  {\cal F} (\beta_j) 
\bigotimes _{0<v, v=0\ \sm (n_j)} {\rm Sym} (\overline{N}_{v,0})\\
\hspace*{30mm} \otimes_{0<v} (\det N_v)^{[{p_j v \over n_j}] +(p-1)v +1} 
 \bigotimes_{0<v, v=0\ \sm (n_j)} (\det N_v)^{-1} ,\\

r_{j*}{\cal F}_{p,j}(X) \simeq {\cal F} (\beta_j) 
\bigotimes _{0<v, v=0\ \sm (n_j)} {\rm Sym} ({N}_{v,0})
\otimes_{0<v} (\det N_v)^{[{p_j v \over n_j}]+(p-1)v  +1},\\

r_{j*}F^1_V \simeq S(V_{0}^{\sR}, \otimes_{0<v} (\det V_v)^{-1}) 
\otimes F^1_V(\beta_j)    
\bigotimes _{0<v, v=0\ \sm (n_j)} \Lambda  (V_{v,0}) \\
\hspace*{30mm}
\otimes_{0<v} (\det \overline{V}_v)^{[{p_j v \over n_j}]+(p-1)v  } ,\\
r_{j*}F^2_V \simeq F^2_V(\beta_j) 
\bigotimes _{0<v, v={n_j \over 2}\ \sm (n_j)} \Lambda  (V_{v,0}) 
\otimes_{0<v} (\det \overline{V}_v)
^{[{p_j v \over n_j} + {1 \over 2}]+(p-1)v } .
\end{array} 
\end{eqnarray}
\end{prop}

{\em Proof} : The proof is similar to that of Proposition 3.1. We divide it 
into several steps.

1) From  (\ref{f1}), (\ref{shif1}), we get 
\begin{eqnarray}\label{shif4}
r_{j*} {\cal F}'_{p-1} (X) = \bigotimes _{\stackrel{0<v,n\in \bN}{0 \leq n \leq (p-1)v} }
{\rm Sym} (N_{v, -n +(p-1)v +{p_j v \over n_j}}) 
\otimes_{0<v} (\det N_v)^{(p-1)v +1} .
\end{eqnarray}

Note that by (\ref{e2}), for $v \in J= \{ v\in \bN|$ There exists $\alpha$ 
such that $N_v\neq 0$ on $F_\alpha \}$, there are no integer in 
$]{p_{j-1} v \over n_{j-1}}, {p_j v \over n_j}[$. 
So for $v\in J$ and $i_0= j-1, \ j$, 
the elements $(v,n)\in \cup_{i=1}^{i_0} I_i^p$ 
are $(v, (p-1)v +1)$, $\cdots, (v, (p-1)v + [{p_{i_0} v \over n_{i_0}}])$. 
Furthermore,
\begin{eqnarray}\label{shif7}\begin{array}{l}
[{p_{j-1} v \over n_{j-1}}]= [ {p_j v \over n_j}] -1
\quad  {\rm if} \quad   v = 0  \quad  \sm (n_j),
\end{array}\end{eqnarray}
\begin{eqnarray}
[{p_{j-1} v \over n_{j-1}}]= [ {p_j v \over n_j}] 
\quad  {\rm if} \quad   v \neq 0  \quad  \sm (n_j). \nonumber
\end{eqnarray}

{}From  (\ref{e3}),  (\ref{shif1}) and  (\ref{shif7}), we have 
\begin{eqnarray}\label{shif6}\begin{array}{l}
r_{j*} 
\Big ( \bigotimes_{(v,n)\in \cup_{i=1}^{j-1} I^p_i} \Big({\rm Sym} ( N_{v,-n}) 
\otimes  \det N_v \Big ) \Big )\\
\hspace*{10mm} = \bigotimes_{(v,n)\in \cup_{i=1}^{j-1} I^p_i} 
{\rm Sym} ( N_{v,-n+(p-1)v + {p_{j} v \over n_{j}}})\otimes _{0<v} 
(\det N_v)^{[{p_{j} v \over n_{j}}] } \\ 
\hspace*{20mm} \otimes_{\stackrel{0<v,}{ v=0\ \sm (n_j)}} (\det N_v)^{-1} ,\\
r_{j*} \Big ( \bigotimes_{(v,n)\in \cup_{i=1}^{j} I^p_i} 
\Big({\rm Sym} ( N_{v,-n}) \otimes  \det N_v \Big )\Big )\\
\hspace*{10mm} = \bigotimes_{(v,n)\in \cup_{i=1}^{j} I^p_i} 
{\rm Sym} ( N_{v,-n+(p-1)v + {p_{j} v \over n_{j}}})\otimes _{0<v} 
(\det N_v)^{[{p_{j} v \over n_{j}}]}. \\
\end{array}\end{eqnarray}

{}From  (3.5), (4.14), (\ref{shif2}), (\ref{shif4}), 
(\ref{shif6}), one easily gets the first two equations in (\ref{shif3}).

2) For $0< n \leq  (p-1)v +{p_j v \over n_j}$, $n\in \bZ$, 
$0\leq i \leq \dim V_v$, (\ref{isomV}) induces a natural  $S^1$-equivariant 
isomorphism of complex vector bundles over $F$,
\begin{eqnarray}\label{shif8}\begin{array}{l}
\Lambda ^i \overline{V}_{v,n-(p-1)v-{p_j v \over n_j} } 
\simeq \Lambda ^{\dim V_v -i} V_{v, -n+(p-1)v+{p_j v \over n_j} } 
\otimes \det  \overline{V}_v.
\end{array}\end{eqnarray}
{}From  (\ref{f2}) and  (\ref{shif8}), as in (\ref{shift1}),
 we obtain the third equality of (\ref{shif3}).

3)  For $0< n \leq (p-1)v +{p_j v \over n_j}+{1 \over 2}$, $n\in \bZ$, 
$0\leq i \leq \dim V_v$, (\ref{isomV}) induces a natural  $S^1$-equivariant
 isomorphism of complex vector bundles over $F$,
\begin{eqnarray}\label{shif9}\begin{array}{l}
\Lambda ^i \overline{V}_{v,n-(p-1)v-{p_j v \over n_j}-{1 \over 2} } 
\simeq \Lambda ^{\dim V_v -i} V_{v, -n+(p-1)v+{p_j v \over n_j}+{1 \over 2} } 
\otimes \det  \overline{V}_v.
\end{array}\end{eqnarray}
{}From (\ref{f2}), (\ref{shif9}), as in (\ref{shift1}),  
we get the last equality of (\ref{shif3}).

The proof of Proposition 4.1 is complete.
\hfill $\blacksquare$

\begin{lemma} Let us introduce the following two line bundles, 
\begin{eqnarray}\label{line5}\begin{array}{l}
L(\beta_j)_1 = L'_1 \otimes_{0<v} (\det N_v )^{[{p_j v \over n_j}]+(p-1)v  +1}
\otimes_{0<v}
(\det \overline{ V}_v)^{[{p_j v \over n_j}]+(p-1)v}    \\
\hspace*{30mm}

\otimes_{0<v, v=0 \ \sm (n_j)} (\det N_v )^{-1} ,\\
L(\beta_j)_2 = L'_2 \otimes_{0<v} (\det N_v)^{[{p_j v \over n_j}] +(p-1)v +1}  
\otimes_{0<v} (\det \overline{ V}_v)
^{[{p_j v \over n_j} + {1 \over 2}]+(p-1)v  } \\
\hspace*{30mm}  \otimes_{0<v, v=0\ \sm (n_j)} (\det N_v )^{-1}.
\end{array}\end{eqnarray}
Then  $L(\beta_j)_1,\ L(\beta_j)_2$ can be extended naturally to  $S^1$-equivariant
complex line bundles which we will still denote by $L(\beta_j)_1,\ L(\beta_j)_2$ respectively over $M(n_j)$.
\end{lemma}

{\em Proof} : We divide the argument into several steps.

 1) Write 
\begin{eqnarray} \label{numb1}
[{p_j v \over n_j}] = {p_j v \over n_j} - {\omega(v)  \over n_j}.
\end{eqnarray}
Note that for $v={n_j \over 2}\ \sm (n_j)$, 
${\omega(v)  \over n_j} = {1 \over 2}$. From (\ref{hyp5}), 
(\ref{line4}) and (\ref{numb1}), we get, formally,
\begin{eqnarray}\label{line6}\begin{array}{l}
L(\beta_j)_1= L^{-(p-1)- p_j /n_j }
 \otimes_{0< v, v\neq {n_j \over 2}\ \sm (n_j)} 
(\det N_v \otimes \det \overline{V}_v)^{-{\omega(v)  \over n_j}} \\
\hspace*{20mm}\otimes_{0<v, {n_j \over 2} <v <n_j\ \sm (n_j) }
(\det N_v \otimes \det \overline{V}_v) \otimes L(n_j)^{{-r(n_j) \over 2n_j}}.
\end{array}\end{eqnarray}
Note that, if $v=v'\ \sm (n_j)$, then $\omega(v) = \omega(v')$. 
Also, for $0< v' <n_j$, $\omega(n_j - v') = n_j - \omega(v')$. 
{}From (\ref{line1}), (\ref{idFV3}), (\ref{idFV4}) and 
(\ref{line6}), formally, we get 
\begin{eqnarray} \qquad\begin{array}{l}
L(\beta_j)_1= L^{-(p-1)-p_j /n_j} \\
\displaystyle{\hspace*{15mm}\otimes 
\Big [ \bigotimes_{0< v< {n_j \over 2}}
 \Big  (\det (N(n_j)_{v}) \otimes \det (\overline{V(n_j)}_{v}) \Big ) 
^{-\omega(v)- r(n_j)v} \Big ]^{1/n_j}.}
\end{array}\end{eqnarray}
We introduce the following line bundle over $M(n_j)$,
\begin{eqnarray}\label{line7}
L^{\omega} (\beta_j)= \bigotimes_{0< v< {n_j \over 2}}
 \Big  (\det (N(n_j)_{v}) \otimes \det (\overline{V(n_j)}_{v}) \Big ) 
^{-\omega(v)- r(n_j)v} .
\end{eqnarray}
Note that $\omega(v) = p_j v \ \sm (n_j)$. If $n_j$ is odd, 
then $r(n_j)=0$, Lemma 4.2
implies that $L^{\omega} (\beta_j)^{1/n_j}$ is well defined  over $M(n_j)$. 
If $n_j$ is even, then  $p_j$ is odd, and $r(n_j)=1$, Lemma 4.2 again implies 
$L^{\omega} (\beta_j)^{1/n_j}$ is well defined over $M(n_j)$. So $L(\beta_j)_1$ extends naturally over $M(n_j)$.

2)   From (\ref{hyp5}),  (\ref{line4}), (\ref{line5}), we get
\begin{eqnarray}\label{line8}
L(\beta_j)_2 = L^{-(p-1)}  \otimes_{0<v} (\det N_v)^{[{p_j v \over n_j}]}  
\otimes_{0<v} (\det \overline{ V}_v)^{[{p_j v \over n_j} + {1 \over 2}]  } 
\otimes_{\stackrel{0<v,}{ {n_j \over 2} < v < n_j \ \sm (n_j)}} 
\det N_v    \nonumber\\
\otimes \Big [ \otimes_{\stackrel{0<v,}{ v={n_j \over 2} \ \sm (n_j)}} 
(\det N_v \otimes \det V_v) \Big ] ^{r(n_j)/2} 
\otimes L(n_j)^{-r(n_j) /2 n_j}. 
\end{eqnarray}

By using the same argument as in 1), one deduces that
\begin{eqnarray}\begin{array}{l}
\displaystyle{\otimes_{0<v} (\det N_v)^{-\omega(v)/n_j} 
\bigotimes_{\stackrel{0<v,}{ {n_j \over 2} <v <n_j\ \sm (n_j)}} \det N_v 
\otimes \Big [ \bigotimes_{ v= {n_j \over 2}\ \sm (n_j)}
 \det N_v \Big ] ^{r(n_j)/2}  }\\
\displaystyle{\hspace*{20mm}= \bigotimes_{0< v < {n_j \over 2} } 
(\det (N(n_j)_{v}) )^{- \omega(v)/n_j}. }
\end{array}\end{eqnarray}
Note that, if there exists $m \in [0,p_j[$ with $ m\in \bN$, such that 
$m \leq p_j v/n_j < (m+ {1\over 2})$, then
\begin{eqnarray}\label{numb2}
[{p_j v \over n_j} + {1 \over 2}] = p_j v/n_j - \omega(v) /n_j.
\end{eqnarray}
While if there exists $m\in ]0,p_j],\ m\in \bN$ such that
 $(m-{1 \over 2}) \leq p_j v/n_j < m$, then
\begin{eqnarray}\label{numb3}
[{p_j v \over n_j} + {1 \over 2}] = p_j v/n_j +(n_j- \omega(v)) /n_j.
\end{eqnarray}

{}From  (\ref{idFV4}), we have the following formal identity,
\begin{eqnarray}\label{line15}
\bigotimes_{\stackrel{0<v'<n_j, m\in \bN} {m < p_j v' /n_j < m+{1\over 2}} }
\bigotimes_{v=v'\ \sm (n_j)} (\det \overline {V}_v)^{-{\omega(v) \over n_j}} 
\bigotimes_{\stackrel{0<v'<n_j, m\in \bN}  
{m-{1 \over 2}<  p_j v' /n_j < m} }
\bigotimes_{v=v'\ \sm (n_j)} (\det \overline {V}_v)
^{{(n_j-\omega(v))\over n_j}}  \nonumber \\
=\bigotimes_{\stackrel{0<v'<n_j/2, m\in \bN} 
{m < p_j v' /n_j < m+{1\over 2}} }
\Big (\bigotimes_{v=v'\ \sm (n_j)} (\det \overline {V}_v)
^{-{\omega(v')\over n_j}}
 \bigotimes_{v=-v'\ \sm (n_j)} (\det \overline {V}_v)
^{{\omega(v') \over n_j}} \Big )\\
\bigotimes_{\stackrel{0<v'<n_j/2, m\in \bN}  
{m-{1 \over 2}<  p_j v' /n_j < m} }
 \left ( \bigotimes_{v=-v' \ \sm (n_j)} (\det \overline {V}_v)
^{-{(n_j-\omega(v'))\over n_j}} 
\bigotimes_{v=v'\ \sm (n_j)} (\det \overline {V}_v)
^{{(n_j-\omega(v'))\over n_j}}  \right )  \nonumber \\
 =\bigotimes_{0<v'<n_j/2} (\det (\overline{V(n_j)}_{v'}) )^{-{\omega(v') \over n_j}} 
\bigotimes_{1\leq m \leq p_j/2} 
\bigotimes_{m-{1 \over 2} < p_j v'/n_j< m} \det (\overline{V(n_j)}_{v'}).\nonumber 
\end{eqnarray}

Recall that for $v\in \bZ$, $p_jv /n_j \in \bZ$ iff $v/n_j \in \bZ$, 
also $p_jv /n_j \in \bZ+{1 \over 2}$ iff $v/n_j \in \bZ+{1 \over 2}$.
{}From (\ref{hyp5}), (\ref{line8})-(\ref{line15}), we get
\begin{eqnarray}\qquad
L(\beta_j)_2 = L^{-(p-1)-{p_j \over n_j}} 
\otimes L^{\omega}(\beta_j)^{{1 \over n_j}} 
\bigotimes_{1\leq m \leq p_j/2} 
\bigotimes_{m-{1 \over 2}<  p_j v' /n_j < m} 
\det (\overline{V(n_j)}_{v'}). 
\end{eqnarray}

The proof of Lemma 4.3 is complete.\hfill $\blacksquare$\\

Let us denote by  
\begin{eqnarray} \begin{array}{l}
\varepsilon_1 = {1 \over 2} \sum_{0<v} (\dim N_v - \dim V_v) \Big [ ( [{p_j v \over n_j}] + (p-1)v) ([{p_j v \over n_j}] + (p-1)v+1)\\
\hspace*{20mm}-( {p_j v \over n_j} + (p-1) v) 
\Big (2 \Big ([{p_j v \over n_j}] + (p-1) v \Big )+1 \Big ) \Big ],\\
\varepsilon_2 = {1 \over 2} \sum_{0<v} (\dim N_v )\Big [
([{p_j v \over n_j}] + (p-1)v) ([{p_j v \over n_j}] + (p-1)v+1)\\
\hspace*{20mm}-( {p_j v \over n_j} + (p-1) v) 
\Big (2([{p_j v \over n_j}] + (p-1) v)+1 \Big ) \Big ]\\
\hspace*{10mm}
-{1 \over 2}\sum_{0<v} (\dim V_v) \Big [ ( [{p_j v \over n_j}+ {1\over 2}]
 + (p-1)v)^2 \\
\hspace*{20mm}
-2 ( {p_j v \over n_j} + (p-1)v)([{p_j v \over n_j}+ {1\over 2}] 
+ (p-1)v) \Big ].
\end{array}\end{eqnarray}
Then  $\varepsilon_1,\ \varepsilon_2$ are locally constant functions on $F$.

Recall that, if $E$ is a $S^1$-equivariant vector bundle over $M$, then 
the weight of the $S^1$ action 
on $\Gamma(F,E)$ is given by the action $\bJ_H$ (cf. Section  3.1).
\begin{prop}  For $i=1,\ 2$, the isomorphisms induced by (4.20) and (\ref{shif3}),
\begin{eqnarray} \label{rel0}\qquad \begin{array}{l}
r_{i1} : S(TY, \otimes_{0<v} (\det N_v)^{-1}) \otimes {\cal F}_{p,j-1}(X) 
\otimes F^i_V \to \\
\hspace*{10mm}S(U_i,L_i)' \otimes {\cal F}(\beta_j) \otimes F^i_V(\beta_j) 
\otimes L(\beta_j)_i
\otimes _{\stackrel{0<v,}{v=0 \sm (n_j)}} {\rm Sym} (\overline{N}_{v,0}),\\
r_{i2} : S(TY, \otimes_{0<v} (\det N_v)^{-1}) \otimes {\cal F}_{p,j}(X) 
\otimes F^i_V \to \\
\hspace*{10mm}S(U_i,L_i)' \otimes {\cal F}(\beta_j) \otimes F^i_V(\beta_j) 
\otimes L(\beta_j)_i
\otimes _{\stackrel{0<v,}{v=0\ \sm (n_j)}}( {\rm Sym} ({N}_{v,0}) 
\otimes \det N_v),
\end{array}\end{eqnarray}
have 
the following properties : 

1) for $i=1,\ 2$, $\gamma=1,\ 2$, one has 
\begin{eqnarray}  \label{rel2} \begin{array}{l}
 \quad r_{i \gamma} ^{-1} \bJ_ H r_{i \gamma} = \bJ_H,\\
 \quad r_{i \gamma} ^{-1} P r_{i \gamma} = P + ({p_j \over n_j} 
+ (p-1)) \bJ_H + \varepsilon_{i \gamma},
\end{array}\end{eqnarray}
where 
\begin{eqnarray}\begin{array}{l}
\varepsilon_{i1} = \varepsilon_i - e(p, \beta_{j-1}, N),\\
\varepsilon_{i2} = \varepsilon_i - e(p, \beta_j, N).
\end{array}\end{eqnarray}

2) Recall that $o(V(n_j)_{n_j \over 2}^\sR)$ was defined in (\ref{vol5}). Denote by   
\begin{eqnarray}\label{rel3} \begin{array}{l}
\mu_1 = -\sum_{0<v} [{p_j v \over n_j}] \dim V_v + \Delta(n_j,N) 
+ \Delta(n_j, V) \quad \sm (2),\\
\mu_2 = - \sum_{0<v} [{p_j v \over n_j}+ {1 \over 2}] \dim V_v 
+ \Delta(n_j,N)+ o(V(n_j)_{n_j \over 2}^\sR) \quad \sm (2),\\
\mu_3 = \Delta(n_j,N)   \quad \sm (2)  .
\end{array}\end{eqnarray}
Then, for $i=1,\ 2$, $\gamma= 1,\ 2$, one has 
\begin{eqnarray}\label{rel4}\begin{array}{l}
r_{i \gamma} ^{-1} \tau_e r_{i \gamma} = (-1)^{\mu_i} \tau_e,\\
r_{i \gamma}^{-1} \tau_s r_{i \gamma}= (-1)^{\mu_3} \tau_s.
\end{array}\end{eqnarray}
\end{prop}

{\em Proof}:  The first equality of (\ref{rel2}) is trivial.

{}From (2.24) and (\ref{shif7}), one has
\begin{eqnarray}\label{rel10}
e(p,\beta_j, N)= e(p, \beta_{j-1}, N) 
+ \sum_{0<v, v=0\ \sm(n_j)}\Big ((p-1)v + {p_j v \over n_j} \Big ) \dim N_v.
\end{eqnarray}

For $i=1, \ 2$, denote by $\varepsilon_i (V)$  the contribution of $\dim V$ to 
$ \varepsilon_i$. From  (\ref{shif8}), 
when acting on $\bigotimes_{\stackrel{0<v, n\in \sN}
{0< n \leq (p-1)v + {p_j v \over n_j}}}$ $ \Lambda^{i_n} \overline{V}_{v,n}$,
 as in (\ref{tranP1}), we get that
\begin{eqnarray}\label{rel11}\begin{array}{l}
 r_{j*}^{-1} P r_{j*} = P + ((p-1)+ {p_j \over n_j}) \bJ_H   \\
\hspace*{20mm}+  \sum _{\stackrel{0<v,n\in \bN,}{ 0< n \leq (p-1)v+ {p_j v \over n_j}}} 
(\dim V_v ) (-n + (p-1)v+ {p_j v \over n_j})  \\
\hspace*{15mm}= P + ((p-1)+ {p_j \over n_j}) \bJ_H + \varepsilon_1 (V) 
- {1 \over 2} ((p-1)+ {p_j \over n_j}) d'(V).
\end{array}\end{eqnarray}

Similarly, by (\ref{shif9}), when acting 
on $\bigotimes_{\stackrel{0<v, n\in \bN}{0< n \leq (p-1)v + 
{p_j v \over n_j} +{1 \over 2}} } \Lambda^{i_n} \overline{V}_{v,n- {1 \over 2}}$, 
as in (\ref{tranP5}),  we get that 
\begin{eqnarray}\label{rel12}\begin{array}{l}
 r_{j*}^{-1} P r_{j*} = P + ((p-1)+ {p_j \over n_j}) \bJ_H \\ 
\hspace*{20mm}+
\sum _{\stackrel{0<v,n\in \bN,}{ 0< n \leq (p-1)v+ {p_j v \over n_j}+{1 \over 2}}} 
(\dim V_v)  (-n + (p-1)v+ {p_j v \over n_j}+{1 \over 2}) \\
\hspace*{15mm}= P + ((p-1)+ {p_j \over n_j}) \bJ_H + \varepsilon_2 (V) .
\end{array}\end{eqnarray}

{}From  (\ref{tranP2}),  (\ref{shif4})-(\ref{shif6}), (\ref{rel10}),
(\ref{rel11}), (\ref{rel12}), and by proceeding as in the proof of 
Proposition 3.2, one deduces easily the second equation of (\ref{rel2}).

3) From  (\ref{shif8}),  on $\bigotimes_{\stackrel{0<v, n\in \sN}{0< n \leq (p-1)v + {p_j v \over n_j}}} \Lambda^{i_n} \overline{V}_{v,n}$, we have
\begin{eqnarray}\label{rel8}
r_{j*}^{-1} \tau_e r_{j*} = (-1)^ {- \Sigma_{0<v}( [{p_j v \over n_j}]+(p-1)v) \dim V_v } \tau_e.
\end{eqnarray}
{}From (\ref{shif9}), on $\bigotimes_{\stackrel{0<v, n\in \sN}{0< n \leq (p-1)v + {p_j v \over n_j}+ {1\over 2}}} \Lambda^{i_n} \overline{V}_{v,n-{1 \over 2}}$, we have
\begin{eqnarray}\label{rel9}
r_{j*}^{-1} \tau_e r_{j*} = 
(-1)^ {- \Sigma_{0<v} ([{p_j v \over n_j}+ {1\over 2}]+(p-1)v) \dim V_v } 
\tau_e.
\end{eqnarray}
As we don't change the $\bZ_2$-grading in the rest part of (\ref{rel0}),
from (\ref{tranP10}), (\ref{idFV11}), (\ref{rel8}), 
 (\ref{rel9}) and the discussion following (\ref{vol5}), 
 we get (\ref{rel4}). 

The proof of Proposition 4.2 is complete.
\hfill $\blacksquare$\\

The following Lemmas 4.4 and 4.6 were essentially proved in 
[{\bf T}, Lemmas 9.6, 9.7].

\begin{lemma} For each connected component $M'$ of $M(n_j)$, 
$\varepsilon_1$, $\varepsilon_2$ are independent on the connected 
component of $F$ in $M'$.
\end{lemma}

{\em Proof} : From  (2.12), (\ref{idFV4}), (\ref{idFV6}), 
 (\ref{numb1}) and (4.45),  we get  
\begin{eqnarray}\label{const1}
\varepsilon_1= {1\over 2} \sum_{0\leq v' < n_j} \sum_{0<v, v = v'\ \sm (n_j)}
(\dim N_v - \dim V_v) \Big [ -( {p_j v \over n_j} + (p-1)v)^2 
- {\omega(v')(n_j -\omega( v')) \over n_j^2} \Big ] \nonumber \\
= (p-1 + {p_j \over n_j})^2 e 
-{1 \over 16} \Big ( \dim_\sR N(n_j)_{n_j\over 2}^{ \sR} 
-\dim_\sR V(n_j)_{n_j\over 2}^{ \sR}  \Big )\\
- {1 \over 2} \sum_{0<v' < n_j/2} 
\Big ( \dim N(n_j)_{v'} -\dim V(n_j)_{v'} \Big )
 { \omega(v')(n_j - \omega(v')) \over n_j^2}.\nonumber 
\end{eqnarray}
By combining (4.11), (4.13), (4.16),  (\ref{numb2}), (\ref{numb3})
and (4.45), as in (\ref{const1}), we get
\begin{eqnarray}\label{const2}
\begin{array}{l}
\displaystyle{
\varepsilon_2 =\varepsilon_1 -{1 \over 2} \sum_{0\leq m <p_j} 
\sum_{m< p_j v'/n_j < m+ {1\over 2} } 
(\dim V(n_j)_{v'} ){\omega(v') \over n_j}   }\\
\displaystyle{\hspace*{10mm}
-{1 \over 2} \sum_{0 <  m \leq p_j} \sum_{m-{1\over 2}< p_j v'/n_j < m } 
(\dim V(n_j)_{v'} ){(n_j- \omega(v'))\over n_j} 
-{1 \over 8} \dim_\sR  V(n_j)_{n_j\over 2}^{ \sR} . }
\end{array}\end{eqnarray}

The proof of Lemma 4.4 is complete.\hfill $\blacksquare$\\

The following Lemma  was proved in [{\bf BT}, Lemma 9.3].

\begin{lemma} Let $M$ be a smooth  manifold on which $S^1$ acts. 
Let $M'$ be a connected  component of $M(n_j)$, 
the fixed point set of $\bZ_{n_j}$ subgroup of $S^1$ on $M$. 
Let $F$ be the fixed point set of $S^1$-action on $M$. 
Let $V\to M$ be a real, oriented, even dimensional vector bundle to 
which the $S^1$-action on $M$ lifts. Assume that $V$ is spin over $M$.
 Let $p_j \in ]0,n_j[$,  $p_j\in \bN$ 
and $(p_j, n_j)= 1$, then
\begin{eqnarray}\label{const3}
\sum_{0<v} (\dim V_v )[{p_j v \over n_j}] + \Delta (n_j, V) \quad \sm (2)
\end{eqnarray}
is independent on the connected components of $F$ in $M'$.
\end{lemma}

Recall that the number $d'(\beta_j,N)$ has been defined in (2.24).

\begin{lemma} For each connected component $M'$ of $M(n_j)$,
$d'(\beta_j, N) + \mu_i \ \sm (2) $ $(i=1,\ 2,\ 3)$ 
is independent on the connected component of $F$ in $M'$.
\end{lemma}

{\em Proof }:   The assertions for the cases $i=1, \ 3$ follow 
immediately from Lemma 4.5 by replacing $V$ with $TX$. 

 From  (4.41) and (4.42), we have 
\begin{eqnarray}\label{const4}
\begin{array}{l}
\displaystyle{
o(V(n_j)_{n_j/2}^{ \sR}) + \sum_{0<v} (\dim V_v) [{p_j v \over n_j} + {1 \over 2}] 
=\sum_{0<v} (\dim V_v )[{p_j v \over n_j}]  + o(V(n_j)_{n_j/2}^{ \sR})  }\\
\displaystyle{ + \sum_{0< m \leq p_j} 
\sum_{\stackrel{0< v' < n_j}{ m -{1 \over 2} < p_j v' /n_j < m}} 
\sum _{\stackrel{0< v,}{ v= v'\ \sm (n_j)}} \dim V_v
+ \sum_{\stackrel{0< v,}{v= {n_j \over 2}\ \sm (n_j) }} \dim V_v.}
\end{array}\end{eqnarray}
By (\ref{idFV6}), the last term in (\ref{const4}), $\dim_\sR V(n_j)_{n_j/2}^{ \sR}$, is a locally constant
 function on $M(n_j)$. 

By (\ref{idFV4}), the third  term on the right side of 
(\ref{const4}) is equal to 
\begin{eqnarray} \label{const5}\qquad \begin{array}{l}
\displaystyle{
\sum_{n_j/2 < v' < n_j} \sum_{0<v, v=v'\ \sm (n_j)} \dim V_v
+ \sum_{0< m \leq p_j} 
\sum_{\stackrel{0< v' < {n_j \over 2}}{ m -{1 \over 2} < p_j v' /n_j < m}} 
\dim V(n_j)_{v'}  \quad \sm (2) }
\end{array}\end{eqnarray}
The last term of (\ref{const5}) is a locally constant function on $M(n_j)$.

By (\ref{vol5}), the first term in (\ref{const5}) is
 $o(V(n_j)_{n_j/2}^{ \sR})$ $ + \Delta(n_j, V) \ \sm (2)$ . 
On the other hand, 
by Lemma 4.5, we know that $\Sigma_{0<v} (\dim V_v) [{p_j v \over n_j}] 
+ \Delta (n_j, V) \
\sm (2) $ is independent on the component of $F$ in $M'$. 

By the above discussion, the left side of (\ref{const4}) is $\sm (2)$
independent on the component of $F$ in $M'$.

The proof of Lemma 4.6 is complete.
\hfill $\blacksquare$

\subsection{ \normalsize Proof of Theorem 2.7}

{}From  (\ref{e6}), (\ref{idFV3}), (\ref{idFV6})  and (\ref{shif7}), we have 
\begin{eqnarray}\label{const8} \qquad
\begin{array}{l}
\displaystyle{
\sum_{0<v} \dim N_v = \sum_{0<v < {n_j \over 2}} \dim N(n_j)_v + {1 \over 2 }
 \dim_\sR N(n_j)_{n_j/2}^{ \sR} + \sum_{0<v, v=0\ \sm (n_j)} \dim N_v,  }\\
\displaystyle{
d'(\beta_j, N) = d'(\beta_{j-1}, N) + \sum_{0<v, v=0\ \sm (n_j)} \dim N_v.}
\end{array} 
\end{eqnarray}

By Lemma 4.6, and  (\ref{const8}), $d'(\beta_{j-1}, N) + \sum_{0<v} \dim N_v 
+ \mu_i \  \sm (2)$ is a constant function on each connected component 
of $M(n_j)$.

{}From  Lemma 4.3, one knows that the Dirac operator 
$D^{X(n_j)} \otimes F(\beta_j) \otimes F^i_V(\beta_j) 
\otimes L(\beta_j)_i $ $(i=1,\ 2)$ is well-defined on $M(n_j)$. Thus,
by using Proposition 4.2, Lemma 4.4,  (\ref{idFV2}), (\ref{idFV12}) and (\ref{const8}), 
for $i=1,\ 2$, $h\in \bZ$, $m\in {1 \over 2} \bZ$,  $\tau= \tau_e$ or $\tau_s$,
and by applying  separately both the first and the second
equations of Theorem 1.2 to each connected component of $M(n_j)$, 
we get the following identity in $K(B)$,
\begin{eqnarray} \label{last}\\
\begin{array}{l}
\sum_\alpha (-1)^{d'(\beta_{j-1}, N) +\sum_{0<v} \dim N_v } 
{\rm Ind}_\tau (D^{Y_\alpha} 
\otimes {\cal F}_{p, j-1}(X)\otimes  F^i_V, 
m +e(p, \beta_{j-1}, N), h)\\
= \sum_\beta (-1)^{d'(\beta_{j-1}, N) + \sum_{0<v} \dim N_v + \mu}
 {\rm Ind}_\tau (D^{X(n_j)} \otimes F(\beta_j) \otimes F^i_V(\beta_j) 
\otimes L(\beta_j)_i, \\
\hspace*{30mm}m+ \varepsilon_i + ({p_j\over n_j} +(p-1)) h, h)\\
=\sum_\alpha (-1)^{d'(\beta_{j}, N)+\sum_{0<v} \dim N_v  } 
{\rm Ind}_\tau (D^{Y_\alpha} 
\otimes  {\cal F}_{p, j}(X)\otimes  F^i_V, 
m+e(p, \beta_{j}, N), h).
\end{array} \nonumber
\end{eqnarray}
Here $\sum_\beta$ means the sum over all the connected components of $M(n_j)$.
 In (\ref{last}), if $\tau=\tau_s$, then $\mu= \mu_3$; if $\tau= \tau_e$, 
then $\mu = \mu_i$.

The proof of Theorem 2.7 is complete.\hfill $\blacksquare$\\

\centerline{------------------------}
\vskip 6mm

Kefeng LIU,
Department of Mathematics, Stanford University, Stanford, CA 94305, USA.

{\em E-mail address}: kefeng@math.stanford.edu

\vskip 6mm

Xiaonan MA,
Humboldt-Universitat zu Berlin, Institut f\"ur Mathematik, unter den Linden 6,
D-10099 Berlin, Germany.

{\em E-mail address}: xiaonan@mathematik.hu-berlin.de

\vskip 6mm

Weiping ZHANG,
Nankai Institute of Mathematics, Nankai university,
Tianjin 300071, P. R. China.

{\em E-mail address}: weiping@nankai.edu.cn

\end{document}